\documentclass[11pt,reqno]{amsart}
\usepackage[letterpaper,margin=1in]{geometry}
\usepackage{amsmath,amssymb,amsthm,mathtools}
\usepackage[hidelinks]{hyperref}
\allowdisplaybreaks[1]
\newtheorem{theorem}{Theorem}[section]
\newtheorem{lemma}[theorem]{Lemma}
\newtheorem{proposition}[theorem]{Proposition}
\newtheorem{corollary}[theorem]{Corollary}
\theoremstyle{definition}
\newtheorem{definition}[theorem]{Definition}
\newtheorem{example}[theorem]{Example}
\newtheorem{problem}[theorem]{Problem}
\theoremstyle{remark}
\newtheorem{remark}[theorem]{Remark}
\numberwithin{equation}{section}

\makeatletter
\@namedef{subjclassname@2020}{\textup{2020} Mathematics Subject Classification}
\makeatother

\title[Kadison's open problem and $R^*$-algebras]
 {Kadison's Open Problem and $R^*$-Algebras}
\author{Bing Yang}
\address{School of Science, Xi'an Polytechnic University, Xi'an 710048, P.\ R.\ China}
\email{yangbmath@foxmail.com}
\thanks{Bing Yang is the corresponding author.}
\author{Chengjun Hou}
\address{School of Mathematical Sciences, Yangzhou University, Yangzhou 225002, P.\ R.\ China}
\email{cjhou@yzu.edu.cn}
\author{Xiaochun Fang}
\address{School of Mathematical Sciences, Tongji University, Shanghai 200092, P.\ R.\ China}
\email{xfang@tongji.edu.cn}

\begin{document}
\begin{abstract}
We study $R^*$-subalgebras of von Neumann algebras. In finite
homogeneous type $\mathrm I_n$, every $R^*$-subalgebra is locally
finite-dimensional. We also give a trace criterion for local finite
dimensionality in the semifinite setting under a finite-trace
hypothesis, and prove an approximate halving property for
projections in strongly dense $R^*$-subalgebras of type
$\mathrm{II}_1$ factors. For three projections in an $R^*$-subalgebra
of a type $\mathrm{II}_1$ factor, the reduced word-length filtration
gives bounds on layer growth, criteria for stabilization, and
uniform bounds on minimal-polynomial degrees at each fixed word
length. Finally, we give a negative answer to Kadison's fifteenth
problem by constructing an infinite-dimensional unital $*$-algebra
generated algebraically by three projections in which every
selfadjoint element has finite spectrum. The resulting $R^*$-algebra
embeds unitally into every type $\mathrm{II}_1$ factor.
\end{abstract}
\maketitle

\section{Introduction and preliminaries}\label{sec:introduction}

At the celebrated 1967 conference in Baton Rouge, Kadison proposed a list of twenty open problems in operator algebras. The list circulated widely and influenced the development of the theory of operator algebras. In 2003, Ge surveyed these problems \cite{Ge}. The fifteenth problem on the list concerns algebraically generated
$*$--algebras of operators on a Hilbert space

\begin{problem}[Kadison, {\cite[Question~15]{Ge}}]
If a selfadjoint operator algebra is finitely generated as an
algebra and every selfadjoint element has finite spectrum, must
the algebra be finite-dimensional?
\end{problem}

Sz\H{u}cs and Tak\'acs \cite{Sz} related this question to Kurosh's
problem and established several equivalent conditions for every
selfadjoint element of a $*$-subalgebra $A$ of $B(H)$ to have finite
spectrum. Mori \cite[Theorem~2.1]{Mor} showed that these conditions
are also equivalent to von Neumann regularity and called the
resulting operator algebras $R^*$-algebras.

\begin{definition}[{\cite[Theorem~2.1]{Mor}}]
Let $A$ be a $*$-subalgebra of $B(H)$. We call $A$ an
$R^*$-algebra if it satisfies any of the following equivalent
conditions.
\begin{enumerate}
 \item $A$ is von Neumann regular, i.e. for every $a\in A$, there is
       $b\in A$ with $aba=a$.
 \item $A$ is algebraic, i.e. For every $a\in A$, the algebra
       $\mathcal A(a)$ generated by $a$ is finite-dimensional.
 \item Every operator in $A$ has closed range.
 \item Every operator in $A$ has finite spectrum.
 \item Every selfadjoint operator in $A$ has closed range.
 \item Every selfadjoint operator in $A$ has finite spectrum.
\end{enumerate}
\end{definition}

In this terminology, Kadison's fifteenth problem asks whether
every $R^*$-algebra is locally finite-dimensional, or equivalently
ultramatricial \cite[Lemma~3.10]{Mor}. Mori \cite[Theorem~4.3]{Mor}
also established a correspondence between commutative
$R^*$-algebras and generalized Boolean algebras. In particular, a
commutative $R^*$-algebra $\mathcal A$ is unital if and only if
its projection lattice $\mathcal P(\mathcal A)$ is a Boolean
algebra. This suggests viewing $R^*$-algebras as noncommutative
analogues of generalized Boolean algebras, just as $C^*$-algebras
and von Neumann algebras are noncommutative analogues of
locally compact spaces and measure spaces, respectively.

An and Gao \cite{An} studied atomic $R^*$-algebras and
$C^*$-uniqueness for pre-$C^*$-algebras and identified several
classes for which Kadison's problem has an affirmative answer.
They also proved that any two projections in an $R^*$-algebra generate
a finite-dimensional algebra \cite[Lemma~2.7]{An} and asked whether
this conclusion extends to arbitrary finite family of projections
\cite[Question~2.1]{An}. Thus, in any counterexample, the smallest possible number of projection generators is three.

Since all projections in an $R^*$-algebra are finite \cite[Section~6]{Mor}, it is natural to
studies $R^*$-subalgebras in semifinite von Neumann algebras.
In Section~\ref{sec:von-neumann}, using polynomial identities, we
prove local finite dimensionality in type $\mathrm I_n$ for
$1 \leq n<\infty$. We also give a trace criterion for local finite
dimensionality when the projections of the subalgebra have finite
trace, and prove that projections in strongly dense
$R^*$-subalgebras of type $\mathrm{II}_1$ factors can be halved up
to a remainder of arbitrarily small trace.

Section~\ref{sec:word-filtration} studies the algebra generated by
three projections in an $R^*$-subalgebra of a type $\mathrm{II}_1$
factor $\mathcal M$ with faithful normal tracial state $\tau$.
In $L^2(\mathcal M,\tau)$, we filter this algebra according to reduced word
length and use $\tau$ to measure the successive filtration layers. We obtain
growth bounds, compression trace estimates, and criteria for
stabilization. We also reduce the finite-dimensionality question
to a corner of the generated algebra.

Section~\ref{sec:counterexample} constructs an infinite-dimensional
$R^*$-algebra generated by three projections. We use the
Clifford generators in the Pauli-matrix form of the Jordan--Wigner
representation \cite{JordanWigner} to label the
edges of triangular bouquets and  then build a layered construction
from the resulting underlying matrices. Schur complements yield
a common finite spectral set across all finite-dimensional
coordinates for every fixed selfadjoint matrix polynomial.
Theorem~\ref{thm:main} gives a negative answer to Kadison's
finite-generation problem. This also shows that
$R^*$-algebras form a significant class of operator algebras.
Corollary~\ref{cor:II1-embedding} shows
that the example embeds unitally into every type $\mathrm{II}_1$
factor. We refer to \cite{Kad2,Tak,Bla} for background on von
Neumann algebras.

\section{\texorpdfstring{$R^*$}{R*}-subalgebras of von Neumann algebras}
\label{sec:von-neumann}

For \(\mathcal F\subseteq B(H)\), let \(\mathcal A(\mathcal F)\) and
\(\mathcal A^*(\mathcal F)\) denote, respectively, the algebra and the \(*\)-algebra
generated by \(\mathcal F\) algebraically. We include the
identity when it is specified among the generators. An algebra is locally
finite-dimensional if each finite subset generates a finite-dimensional
subalgebra. For a \(*\)-algebra, one can equivalently require each finite
subset to generate a finite-dimensional \(*\)-subalgebra, since adjoining adjoints of a finite set again yields a finite set.

\begin{lemma}[Amitsur--Levitzki, \cite{Dre}]\label{lem:1}
If \(Z\) is a commutative ring, the standard polynomial
\[
s_{2n}(x_1,\ldots,x_{2n})
=\sum_{\pi\in\mathfrak S_{2n}}\operatorname{sgn}(\pi)
x_{\pi(1)}\cdots x_{\pi(2n)}
\]
vanishes identically on \(M_n(Z)\).
\end{lemma}

\begin{lemma}[Kaplansky, {\cite[Theorem~6.1]{Kap}}]\label{lem:2}
An algebraic algebra satisfying a polynomial identity is locally
finite-dimensional.
\end{lemma}

\begin{theorem}\label{thm:finite-type-I}
Let \(\mathcal M\) be a von Neumann algebra of type \(\mathrm I_n\) with
\(1\leq n<\infty\). Then every \(R^*\)-subalgebra of \(\mathcal M\) is
locally finite-dimensional.
\end{theorem}

\begin{proof}
The type decomposition gives
\[
\mathcal M\cong \mathcal A_n\,\overline{\otimes}\,B(H_n)
\cong M_n(\mathcal A_n),
\qquad \dim H_n=n,
\]
where \(\mathcal A_n\) is an abelian von Neumann algebra and \(H_n\) is
an \(n\)-dimensional Hilbert space. By Lemma~\ref{lem:1}, the standard
identity \(s_{2n}=0\) holds on \(\mathcal M\) and hence on its subalgebras.
Every \(R^*\)-algebra is algebraic by \cite[Theorem~2.1]{Mor}, so the
conclusion follows from Lemma~\ref{lem:2}.
\end{proof}

The two-projection result below is due to An and Gao
\cite[Lemma~2.7]{An}. We give a proof using a central element.

\begin{theorem}\label{thm: 2.4}\label{thm:3}
If \(p,q\) are projections in an \(R^*\)-algebra, then
\(\mathcal A^*(1,p,q)\) is finite-dimensional.
\end{theorem}

\begin{proof}
Put \(z=(p-q)^2=p-pq-qp+q\). Then \(z^*=z\), and direct multiplication gives
\[
pz=zp=p-pqp,\qquad qz=zq=q-qpq.
\]
Let \(z=\sum_{i=1}^r\lambda_i e_i\) be its finite spectral resolution,
with \(\sum_i e_i=1\). Then each $e_i$ lies in the center of $A(1,p,q)$. Now we consider the subalgebra $e_iA(1,p,q)e_i$. Obviously,
\[
\begin{aligned}
pqe_i &= pe_i+qe_i-qpe_i-\lambda_i e_i,\\
pqpe_i &= (p-pzp)e_i=p(1-z)e_i=(1-\lambda_i)pe_i,\\
qpqe_i &= (1-\lambda_i)qe_i.
\end{aligned}
\]
Therefore,
every operator in $e_iA(1,p,q)e_i$ can be expressed as a linear combination of the elements of set of $\{e_i, pe_i, qe_i, qpe_i\}$. Thus $\dim A(1,p,q)\leq 4r<\infty $.
\end{proof}

Following \cite[Definition~2.2]{An}, we call an \(R^*\)-algebra
\emph{atomic} if every nonzero projection in it dominates a minimal
nonzero projection. Mori's \emph{purely atomic}
condition requires every projection to be a finite sum of atoms and
implies local finite-dimensionality \cite[Lemma~5.2]{Mor}.

\begin{theorem}\label{thm:atomic-enlargement}
Every represented \(R^*\)-algebra is contained in an atomic
\(R^*\)-algebra on the same Hilbert space. Consequently, all atomic
\(R^*\)-algebras are locally finite-dimensional if and only if all
\(R^*\)-algebras are locally finite-dimensional.
\end{theorem}

\begin{proof}
Let \(\mathcal B\subseteq B(H)\) be an \(R^*\)-algebra, let
\(\mathcal F(H)\) be the ideal of finite-rank operators, and put
\(\mathcal A=\mathcal B+\mathcal F(H)\). This is a \(*\)-subalgebra.
Every selfadjoint \(a\in\mathcal A\) has a decomposition \(a=b+f\)
with \(b\in\mathcal B\) and \(f\in\mathcal F(H)\) both selfadjoint, by
taking the selfadjoint parts of a decomposition. Write
\(b=\sum_{j=1}^r\lambda_j e_j\) with \(\sum_j e_j=I_H\), and put
\[
E=\sum_{j=1}^r e_j(\operatorname{ran} f).
\]
This is a finite-dimensional subspace invariant under the projections
\(e_j\), so it reduces \(b\). It contains \(\operatorname{ran} f\).
Selfadjointness of \(f\) then gives \(f|_{E^\perp}=0\), and \(E\) also
reduces \(f\). Thus \(a=a|_E\oplus b|_{E^\perp}\) has finite spectrum.
This proves that \(\mathcal A\) is an \(R^*\)-algebra. If \(p\in\mathcal A\)
is a nonzero projection, the projection onto any line in \(pH\)
is a rank-one element of \(\mathcal A\) below \(p\), so \(\mathcal A\)
is atomic.

If all atomic \(R^*\)-algebras are locally finite-dimensional,
then so are \(\mathcal A\) and its subalgebra \(\mathcal B\).
The converse follows by restricting the assertion to
atomic \(R^*\)-algebras.
\end{proof}

\begin{theorem}\label{thm:trace-lower-bound}
Let \(\mathcal M\) be a finite von Neumann algebra with a faithful normal
tracial state \(\tau\), and let \(\mathcal A\subseteq\mathcal M\) be an
\(R^*\)-subalgebra. Then \(\mathcal A\) is locally finite-dimensional if
and only if, for every finite \(\mathcal F\subseteq\mathcal A\), there is
\(\delta_{\mathcal F}>0\) such that
\(\tau(p)\geq\delta_{\mathcal F}\) for every nonzero projection
\(p\in\mathcal A^*(\mathcal F)\).
\end{theorem}

\begin{proof}
Fix a finite subset \(\mathcal F\) and put \(\mathcal B=\mathcal A^*(\mathcal F)\).
if \(\mathcal B=\{0\}\), the conclusion holds trivially, so we assume
\(\mathcal B\ne\{0\}\). If \(\mathcal A\) is locally finite-dimensional, then
\(\mathcal B\) is a finite-dimensional \(C^*\)-algebra, and
\[
\mathcal B\cong\bigoplus_{j=1}^r M_{n_j}(\mathbb C).
\]
In this decomposition the faithful trace has the form
\[
\tau(b)=\sum_{j=1}^r c_j\operatorname{Tr}_{n_j}(b_j),
\qquad b=(b_1,\ldots,b_r)\in\mathcal B,\qquad c_j>0,
\]
where \(\operatorname{Tr}_{n_j}\) denotes the unnormalized matrix trace.
Set \(\delta_{\mathcal F}=\min_{1\leq j\leq r}c_j\).
A nonzero projection \(p\in\mathcal B\) has a nonzero matrix block and
therefore satisfies \(\tau(p)\geq\delta_{\mathcal F}\).

Conversely, assume the lower bound \(\delta_{\mathcal F}\). Let
\(y=\sum_{x\in\mathcal F}(x^*x+xx^*)\). Denote by
\(\sigma(y)\setminus\{0\}:=\{\lambda_1,\ldots,\lambda_r\}\),
and define the Lagrange interpolation polynomial by \(p(t)=1-\prod_{i=1}^{r}(1-t/\lambda_i)\) on
\(\sigma(y)\cup\{0\}\). Then
the support projection of $y$ is \(e_{\mathcal F}=p(y)\). Thus \(e_{\mathcal F}\in\mathcal B\), even if
\(\mathcal B\) is generated without adjoining \(1\).
From \(x^*x\leq y\) and \(xx^*\leq y\) we get
\(x(1-e_{\mathcal F})=x^*(1-e_{\mathcal F})=0\) for \(x\in\mathcal F\).
Hence \(e_{\mathcal F}x=xe_{\mathcal F}=x\), so \(e_{\mathcal F}\) is
the identity of \(\mathcal B\).

If \(p_1,\ldots,p_n\) are pairwise orthogonal nonzero projections
in \(\mathcal B\), then
\[
n\delta_{\mathcal F}
\leq \sum_{i=1}^n\tau(p_i)
=\tau\!\left(\sum_{i=1}^n p_i\right)
\leq \tau(e_{\mathcal F})\leq 1,
\]
so \(n\leq1/\delta_{\mathcal F}\). In particular, any orthogonal family
\((q_i)_{i\in I}\) of nonzero projections in \(\mathcal B\) has finite
cardinality. Otherwise, a finite subset \(I_0\subseteq I\) with
\(|I_0|>1/\delta_{\mathcal F}\) would give
\[
1\geq \tau\!\left(\sum_{i\in I_0}q_i\right)
=\sum_{i\in I_0}\tau(q_i)
\geq |I_0|\delta_{\mathcal F}>1.
\]
Choose such a family \(q_1,\ldots,q_m\) of maximal cardinality.
Its sum \(e=\sum_{i=1}^m q_i\) must equals \(s(y)\);otherwise
\(s(y)-e\) could be added. Each \(q_i\) is minimal, since a proper
nonzero subprojection would allow one member to be split into two.
For each \(i\), we have
\(q_i\mathcal Bq_i=\mathbb C q_i\).
If \(a,b\in q_i\mathcal Bq_j\) with \(\|a\|=\|b\|=1\), then
\[
a^*a,b^*b\in\mathbb Cq_j,\qquad
\|a^*a\|=\|b^*b\|=1,\qquad a^*a=b^*b=q_j.
\]
Also \(ab^*\in q_i\mathcal Bq_i\), so \(ab^*=\lambda q_i\) for some
\(\lambda\in\mathbb C\). Multiplication by \(b\) gives \(a=\lambda b\).
Thus \(\dim(q_i\mathcal Bq_j)\leq1\). As
\(\mathcal B=\sum_{i,j=1}^m q_i\mathcal Bq_j\), we conclude that
\[
\dim\mathcal B\leq m^2\leq \frac{1}{\delta_{\mathcal F}^2}.
\]
This proves local finite-dimensionality of \(\mathcal A\).
\end{proof}

\begin{corollary}\label{cor:semifinite-trace}
Let \(\mathcal M\) have a faithful normal semifinite trace \(\tau\).
Suppose that \(\mathcal A\subseteq\mathcal M\) is an \(R^*\)-subalgebra
and that \(\tau(p)<\infty\) for every projection \(p\in\mathcal A\).
Then the trace criterion of Theorem~\ref{thm:trace-lower-bound}
holds for \(\mathcal A\).
\end{corollary}

\begin{proof}
For a finite \(\mathcal F\) with \(\mathcal A^*(\mathcal F)\ne0\), the
support construction in the proof of
Theorem~\ref{thm:trace-lower-bound} gives an identity projection
\(e\in\mathcal A^*(\mathcal F)\). It has \(0<\tau(e)<\infty\), and
\(e\mathcal M e\) has the faithful normal tracial state
\(\tau(\,\cdot\,)/\tau(e)\). If \(\mathcal A^*(\mathcal F)\) is
finite-dimensional, the trace bound follows by applying
Theorem~\ref{thm:trace-lower-bound} in \(e\mathcal M e\). Conversely,
the assertion follows from Theorem ~\ref{thm:trace-lower-bound}.
\end{proof}

For a von Neumann algebra \(\mathcal M\), let \(\mathcal P(\mathcal M)\)
denote the set of all projections in \(\mathcal M\). If \(\tau\) is a faithful normal tracial
state, put \(\|x\|_2=\tau(x^*x)^{1/2}\). A von Neumann algebra is
continuous if and only if each projection \(p\) can be written
as \(p=q+r\) with \(qr=0\) and \(q\sim r\)
\cite[Lemma~6.5.6]{Kad2}. We prove an approximate form of this
decomposition for strongly dense \(R^*\)-subalgebras of
type \(\mathrm{II}_1\) factors.

\begin{lemma}\label{thm:4}
Let \(\mathcal M\) be a type \(\mathrm{II}_1\) factor with a faithful
normal tracial state \(\tau\), and let \(\mathcal A\) be a strongly
dense \(R^*\)-subalgebra. For every projection \(p\in\mathcal M\) and
every \(\eta>0\), there is a projection \(q\in\mathcal A\) such that
\(\|q-p\|_2<\eta\), where \(\|x\|_2=\tau(x^*x)^{1/2}\).
\end{lemma}

\begin{proof}
By Kaplansky's density theorem, the selfadjoint unit ball of
\(\overline{\mathcal A}^{\|\cdot\|}\) is strongly dense in the
selfadjoint unit ball of \(\mathcal M\).

Since the topology generated by \(\|\cdot\|_2\) is equivalent to the
strong operator topology in the unit ball \({\mathcal M}_1\) of \(\mathcal M\).
Hence there exists a selfadjoint \(a\in\mathcal A\) so that \(\|a-p\|_2<\eta/2\).

Let \(\sigma(a)\setminus\{0\}=\{\mu_1,\ldots,\mu_s\}\), and let $f$ be the Lagrange interpolation polynomial on \(\sigma(a)\cup\{0\}\) given by
\[
f(t)=
\sum_{\substack{1\le i\le s\\\mu_i>1/2}}
\frac{t}{\mu_i}
\prod_{\substack{1\le j\le s\\j\ne i}}
\frac{t-\mu_j}{\mu_i-\mu_j}.
\]
Set \(q=\mathbf 1_{(1/2,\infty)}(a)\).
By the functional calculus,
\((a-1)^2\geq(1-q)/4\) and \(a^2\geq q/4\). Therefore
\[
\begin{aligned}
\|a-p\|_2^2
 &=\tau\bigl(p(a-1)^2\bigr)+\tau\bigl((1-p)a^2\bigr)\\
 &\geq \tfrac14\bigl(\tau(p(1-q))+\tau((1-p)q)\bigr)
 =\tfrac14\|p-q\|_2^2 .
\end{aligned}
\]
Thus \(\|p-q\|_2\leq2\|a-p\|_2<\eta\).
\end{proof}

\begin{theorem}\label{thm:approximate-halving}
Let \(\mathcal M\) be a type \(\mathrm{II}_1\) factor with a faithful
normal tracial state \(\tau\), and let \(\mathcal A\subseteq\mathcal M\)
be a strongly dense \(R^*\)-subalgebra. For every nonzero projection
\(p\in\mathcal A\) and every \(\varepsilon>0\), there are pairwise
orthogonal projections \(p_1,p_2,r\in\mathcal A\) such that
\[
p=p_1+p_2+r,\qquad
p_1\sim_{\mathcal A}p_2,\qquad
\tau(r)<\varepsilon.
\]
Here the equivalence is implemented by a partial isometry in
\(\mathcal A\), and \(p_1,p_2\) can both be taken nonzero.
\end{theorem}

\begin{proof}
In \(p\mathcal M p\), choose projections \(e_0,f_0\) with
\(p=e_0+f_0\), \(e_0f_0=0\), and \(e_0\sim f_0\).
Take a partial isometry \(v\in p\mathcal M p\) with
\(v^*v=e_0\) and \(vv^*=f_0\). Put \(t=\tau(p)/2>0\), and choose
\(\delta>0\) with \(\delta^2<\min\{t,\varepsilon/2\}\).
Since \(p\mathcal A p\) is strongly dense in \(p\mathcal M p\),
 Lemma~\ref{thm:4} implies that there is
a projection \(e\in p\mathcal A p\) such that
\(\|e-e_0\|_2<\delta/4\). Similarly, we can choose \(a\in\mathcal A\) such that \(\|a-v\|_2<\delta/2\).
Put \(x=(p-e)ae\in\mathcal A\). Since \(v=(p-e_0)ve_0\) and
\(\|v\|=1\), we have
\[
\begin{aligned}
x-v
 &= (p-e)ae-(p-e_0)ve_0\\
 &= (p-e)(a-v)e+(p-e)ve-(p-e)ve_0\\
 &\quad+\bigl((p-e)-(p-e_0)\bigr)ve_0\\
 &= (p-e)(a-v)e+(p-e)v(e-e_0)+(e_0-e)ve_0.
\end{aligned}
\]
Consequently, \(\|x-v\|_2\leq\|a-v\|_2+2\|e-e_0\|_2<\delta\).

By \cite[Proposition~2.6]{Mor}, write \(x=u|x|\) with \(u\in\mathcal A\). Set \(p_1=u^*u\) and \(p_2=uu^*\).
Then \(p_1\leq e\) and \(p_2\leq p-e\), and \(u\) implements their
equivalence; hence \(\tau(p_1)=\tau(p_2)\). From \(x(p-p_1)=0\)
and \(v^*v=e_0\) we obtain
\[
\begin{aligned}
\|v-x\|_2^2
 &\geq \|(v-x)(p-p_1)\|_2^2
 =\|v(p-p_1)\|_2^2\\
 &=\tau\bigl((p-p_1)e_0(p-p_1)\bigr)
 =t-\tau(e_0p_1)
 \geq t-\tau(p_1).
\end{aligned}
\]
Thus \(\tau(p_1)\geq t-\|v-x\|_2^2>t-\delta^2>0\), so \(p_1,p_2\)
are both nonzero. With \(r=p-p_1-p_2\), we have
\[
0\leq\tau(r)=2t-2\tau(p_1)
\leq2\|v-x\|_2^2<2\delta^2<\varepsilon.
\]
\end{proof}

The next two examples use left regular representations to obtain
\(R^*\)-algebras whose von Neumann closures are of type
\(\mathrm{II}_1\) and \(\mathrm{II}_\infty\), respectively.

\begin{example}\label{ex:finitary-symmetric}
Let \(\Omega\) be a countably infinite set. For
\(g\in\operatorname{Sym}(\Omega)\), put
\(\operatorname{supp}(g)=\{\omega\in\Omega:g(\omega)\ne\omega\}\).
The finitary symmetric group \(G=\operatorname{Fin}(\Omega)\) consists
of those \(g\) with finite support. A group is locally finite if every
finite subset generates a finite subgroup. Each finite set of elements
of \(G\) acts on the finite union of their supports, so \(G\) is
locally finite. The group algebra of a discrete group is an
\(R^*\)-algebra in the left regular representation exactly when the
group is locally finite \cite[Example~3.9]{Mor}. In particular,
\(\mathbb C[G]\) is an \(R^*\)-algebra. For each finite subgroup
\(H\leq G\), the operators \(\{\lambda_g:g\in H\}\) generate the
finite-dimensional \(C^*\)-subalgebra \(\mathbb C[H]\). Since a finite
set of elements of \(\mathbb C[G]\) has a common finite group support,
\(\mathbb C[G]\) is locally finite-dimensional.

Let \(g\in G\setminus\{1\}\) and put
\(S=\operatorname{supp}(g)
   =\{\omega\in\Omega:g(\omega)\ne\omega\}.\)
Choose pairwise disjoint subsets \(S_1,S_2,\ldots\) of
\(\Omega\setminus S\) with \(|S_n|=|S|\). For each \(n\), fix a
bijection \(\varphi_n:S\to S_n\) and set
\[
 h_n=\prod_{\omega\in S}(\omega\ \varphi_n(\omega))\in G.
\]
The transpositions in this product are disjoint, and \(h_n(S)=S_n\).
Consequently,
\[
 \operatorname{supp}(h_ngh_n^{-1})
 =h_n\bigl(\operatorname{supp}(g)\bigr)
 =S_n.
\]
Since the sets \(S_n\) are nonempty and pairwise disjoint, the conjugates
\(h_ngh_n^{-1}\) are all distinct. Thus \(G\) is ICC, and by \cite[Theorem~6.7.5]{Kad2}
\(L(G)=\mathcal L_{\operatorname{Fin}(\Omega)}\) is a type \(\mathrm{II}_1\) factor.

Choose a sequence of finite subsets of \(\Omega\) such that
\[
 \Omega_1\subseteq\Omega_2\subseteq\cdots,\qquad
 |\Omega_n|<\infty,\qquad
 \Omega=\bigcup_{n\geq1}\Omega_n,
\]
and define
\(
 G_n=\{g\in G:\operatorname{supp}(g)\subseteq\Omega_n\}
 \cong\operatorname{Sym}(\Omega_n).
\)
Every element of \(G\) has finite support, so
\[
 G_1\subseteq G_2\subseteq\cdots,\qquad
 G=\bigcup_{n\geq1}G_n.
\]
In the left regular representation,
\(\mathbb C[G_n]=\operatorname{span}\{\lambda_g:g\in G_n\}\) and \(
 \dim_{\mathbb C}\mathbb C[G_n]=|G_n|<\infty.\)
Therefore
\(L(G)=\overline{\mathbb C[G]}^{\,\mathrm{SOT}}
     =\overline{\bigcup_{n\geq1}\mathbb C[G_n]}^{\,\mathrm{SOT}}.\)
Thus \(L(G)\) is hyperfinite.

For a projection \(p=\sum_{s\in F}a_s\lambda_s\) in
\(\mathbb C[G]\), take \(F\) to be its finite support. This
representation is unique, since applying \(p\) to \(\delta_1\)
recovers its coefficients. With respect to the orthonormal basis
\(\{\delta_g:g\in G\}\), the entries of \(p\) are
\[
\begin{aligned}
p_{s,t}
 &= \left\langle\left(\sum_{l\in F}a_l\lambda_l\right)
    \delta_t,\delta_s\right\rangle
 = \sum_{l\in F}a_l\langle\lambda_l\delta_t,\delta_s\rangle\\
 &=
 \begin{cases}
 a_{st^{-1}},& st^{-1}\in F,\\
 0,& st^{-1}\notin F,
 \end{cases}
 \qquad s,t\in G.
\end{aligned}
\]
Thus the entries are constant on each group diagonal
\(D_l=\{(s,t)\in G\times G:st^{-1}=l\}\), and precisely the
\(|F|\) diagonals indexed by \(l\in F\) have nonzero entries.

Let \(H_F=\langle F\rangle\), so \(H_F\) is a finite subgroup and
\(p\in\mathbb C[H_F]\). Choose a set \(T\subseteq G\) of
representatives for the right cosets of \(H_F\). Then
\[
G=\bigcup_{t\in T}H_Ft,\qquad
\ell^2(G)=\bigoplus_{t\in T}\ell^2(H_Ft).
\]
For \(l,h\in H_F\) and \(t\in T\),
\(\lambda_l\delta_{ht}=\delta_{lht}\in\ell^2(H_Ft)\).
This also holds for the adjoint \(\lambda_l^*=\lambda_{l^{-1}}\),
so each subspace \(\ell^2(H_Ft)\) reduces \(p\). Let
\(\lambda^{H_F}\) be the left regular representation of \(H_F\)
on \(\ell^2(H_F)\) and put
\[
P_{H_F}=\sum_{l\in F}a_l\lambda^{H_F}_l.
\]
For \(t\in T\), the unitary
\(U_t\colon\ell^2(H_F)\to\ell^2(H_Ft)\) defined by
\(U_t\delta_h=\delta_{ht}\) satisfies
\[
U_t^*\bigl(p|_{\ell^2(H_Ft)}\bigr)U_t=P_{H_F}.
\]
Therefore \(p\cong\bigoplus_{t\in T}P_{H_F}\), a direct sum
of copies of the same finite matrix.
\end{example}

\begin{example}\label{ex:semifinite-group}
Let \(G\) be a countable discrete locally finite group and put
\(\mathcal K=\ell^2(G)\). Represent \(\mathbb C[G]\) on \(\mathcal K\)
by the left regular representation, and let \(\mathcal F(\mathcal K)\)
denote the algebra of finite-rank operators on \(\mathcal K\).
Then
\[
\mathcal A=\mathbb C[G]\odot\mathcal F(\mathcal K)
\subseteq B(\mathcal K\otimes\mathcal K)
\]
is a locally finite-dimensional \(R^*\)-algebra, with
\[
\overline{\mathcal A}^{\,\mathrm{WOT}}
=L(G)\,\overline{\otimes}\,B(\mathcal K).
\]
If \(G\) is infinite and ICC, the closure is a type
\(\mathrm{II}_\infty\) factor.
\end{example}

\begin{proof}
Take a finite subset of \(\mathcal A\) and express its elements as
finite sums of elementary tensors \(a_i\otimes b_i\). The union of
the group supports of the \(a_i\) generates a finite subgroup
\(H\leq G\), so every \(a_i\) belongs to \(\mathbb C[H]\). Set
\[
V=\sum_i\bigl(\operatorname{ran} b_i+\operatorname{ran} b_i^*\bigr).
\]
The space \(V\) is finite-dimensional and \(b_i=P_Vb_iP_V\) for
each \(i\), where \(P_V\) is the orthogonal projection onto \(V\).
Thus the chosen finite subset lies in \(\mathbb C[H]\otimes B(V)\),
where \(B(V)\) acts by zero on \(V^\perp\). This is a
finite-dimensional \(C^*\)-subalgebra of \(\mathcal A\), and hence
\(\mathcal A\) is locally finite-dimensional.
Hence \(\mathcal A\) is locally finite-dimensional. For each
\(a\in\mathcal A\), choose a finite-dimensional \(C^*\)-subalgebra
\(\mathcal D\subseteq\mathcal A\) containing \(a\), and let \(e\) be its
identity. Since \(a=eae\), it acts as zero on the orthogonal complement
of the range of \(e\). Thus
\[
 \sigma_{B(\ell^2(G)\otimes\mathcal K)}(a)
 \subseteq \sigma_{\mathcal D}(a)\cup\{0\}.
\]
The spectrum \(\sigma_{\mathcal D}(a)\) is finite because \(\mathcal D\)
is finite-dimensional. Therefore every element of \(\mathcal A\) has
finite spectrum, and \(\mathcal A\) is an \(R^*\)-algebra by
\cite[Theorem~2.1]{Mor}.

Let \(\mathcal N=\overline{\mathcal A}^{\,\mathrm{WOT}}\).
Choose finite-rank orthogonal projections \(P_n\) on \(\mathcal K\)
converging strongly to \(I_{\mathcal K}\). For \(a\in\mathbb C[G]\),
\[
a\otimes P_n \longrightarrow a\otimes I_{\mathcal K}
\qquad\text{strongly}.
\]
In particular, \(\mathcal N\) contains the identity, so it is a
von Neumann algebra. The displayed convergence implies
\(L(G)\otimes I_{\mathcal K}\subseteq\mathcal N\). Moreover,
\(I_{\mathcal K}\otimes\mathcal F(\mathcal K)\subseteq\mathcal A\),
and finite-rank operators are strongly dense in \(B(\mathcal K)\).
Thus \(I_{\mathcal K}\otimes B(\mathcal K)\subseteq\mathcal N\).
These two subalgebras generate
\(L(G)\,\overline{\otimes}\,B(\mathcal K)\), which is weak operator
closed and contains \(\mathcal A\). This proves the stated equality.

If \(G\) is infinite and ICC, then \(L(G)\) is a type
\(\mathrm{II}_1\) factor. Since \(\mathcal K\) is infinite-dimensional,
\(L(G)\,\overline{\otimes}\,B(\mathcal K)\) is a type
\(\mathrm{II}_\infty\) factor.
\end{proof}

\section{Subalgebras generated by three projections in \texorpdfstring{$R^*$}{R*}-algebras of \texorpdfstring{$\mathrm{II}_1$}{II1} factors}
\label{sec:word-filtration}

Let $\mathcal M$ be a type $\mathrm{II}_1$ factor with a faithful
normal tracial state $\tau$, and let $\mathcal A\subseteq\mathcal M$
be an $R^*$-algebra. The inner product
$\langle x,y\rangle=\tau(y^*x)$ on $\mathcal M$ defines the norm
$\|\cdot\|_2$, and $L^2(\mathcal M,\tau)$ denotes the corresponding
Hilbert space completion. The $R^*$-subalgebra generated by two
projections is finite-dimensional by \cite[Lemma~2.7]{An} and
Theorem~\ref{thm: 2.4}. We study the three-projection case using the
word-length filtration.

Fix $p_1,p_2,p_3\in\mathcal P(\mathcal A)$ and put
$\mathcal B=\operatorname{alg}(1,p_1,p_2,p_3)$. If $1\notin\mathcal A$,
we work in the $R^*$-algebra $\mathcal A+\mathbb C1$.
Define
\[
 V_n=\operatorname{span}_{\mathbb C}
 \{p_{i_1}\cdots p_{i_k}:0\leq k\leq n,\ i_j\in\{1,2,3\}\},
 \qquad V_0=\mathbb C1,\quad V_{-1}=\{0\}.
\]
A product of linear subspaces denotes the span of their products.
Splitting words at an appropriate position gives $V_aV_b=V_{a+b}$
for $a,b\geq0$. Let $P_n$ be the orthogonal projection of
$L^2(\mathcal M,\tau)$ onto $V_n$, with $P_{-1}=0$. Set
\[
 Q_n=P_n-P_{n-1},\qquad H_n=V_n\ominus V_{n-1},
 \qquad \mathcal B=\bigcup_{n\geq0} V_n.
\]
Let $\mathcal W_{\leq n}$ denote the formal reduced words of length
at most $n$, including the empty word. A word is reduced if its
adjacent letters are distinct. The cardinality of $\mathcal W_{\leq n}$
is
\[
 u_n=1+\sum_{j=1}^n 3\cdot 2^{j-1}=3\cdot 2^n-2.
\]

\begin{theorem}\label{the:2.1}
For $w\in\mathcal W_{\leq n}$, write $\widehat w=w(p_1,p_2,p_3)$.
Define $\pi_n\colon\mathbb C^{\mathcal W_{\leq n}}\to V_n$ by
$\pi_n((c_w)_w)=\sum_w c_w\widehat w$, and set
\[
 K_n=\ker\pi_n,\qquad
 G_n=\bigl(\tau(\widehat v^*\widehat w)\bigr)_{v,w\in\mathcal W_{\leq n}}.
\]
Then the following hold.
\begin{enumerate}
 \item The matrix $G_n$ is positive semidefinite, and
 $K_n=\ker G_n$ and $\dim V_n=\operatorname{rank}G_n$.
 \item With $h_{n+1}=\dim V_{n+1}-\dim V_n$ and
 $k_{n+1}=\dim K_{n+1}-\dim K_n$, one has
 $h_{n+1}+k_{n+1}=3\cdot2^n$.
 \item Write the evaluations of the reduced words of length exactly
 $n+1$ as $a_1,\ldots,a_{3\cdot2^n}$, indexed by the formal words.
 Set $r_j=(I-P_n)a_j$ and $S_n=(\tau(r_i^*r_j))_{ij}$.
 Then $h_{n+1}=\operatorname{rank}S_n$ and
 $k_{n+1}=\dim\ker S_n$.
\end{enumerate}
\end{theorem}

\begin{proof}
For $c=(c_w)_w$ one has
\[
 c^*G_nc
 =\tau\!\left(\left(\sum_wc_w\widehat w\right)^*
                   \left(\sum_wc_w\widehat w\right)\right)
 =\|\pi_n(c)\|_2^2\geq0.
\]
Thus $G_n$ is positive semidefinite, and faithfulness of $\tau$
gives $K_n=\ker G_n$. The identity $p_i^2=p_i$ allows every word
to be reduced without increasing its length, so $\pi_n$ is onto.
Rank-nullity gives
\[
 \dim V_n=u_n-\dim K_n=u_n-\dim\ker G_n=\operatorname{rank}G_n.
\]
Consequently,
\[
 \begin{aligned}
 h_{n+1}+k_{n+1}
 &=\bigl(\dim V_{n+1}+\dim K_{n+1}\bigr)
   -\bigl(\dim V_n+\dim K_n\bigr)\\
 &=u_{n+1}-u_n=3\cdot2^n.
 \end{aligned}
\]
The residuals $r_j$ span $H_{n+1}$ since the words of length at
most $n$ is not in $H_{n+1}$. Applying the Gram matrix argument to
the map $(c_j)_j\mapsto\sum_jc_jr_j$ gives
$\operatorname{rank}S_n=\dim H_{n+1}=h_{n+1}$. Thus
$\dim\ker S_n=3\cdot2^n-h_{n+1}=k_{n+1}$.
\end{proof}

\begin{corollary}\label{cor:word-stabilization}
For every $N\geq0$,
\[
 S_N=0
 \iff\operatorname{rank}G_{N+1}=\operatorname{rank}G_N
 \iff V_{N+1}=V_N
 \iff \mathcal B=V_N.
\]
Such an $N$ exists if and only if $\mathcal B$ is finite-dimensional.
\end{corollary}

\begin{proof}
The first three conditions are equivalent by Theorem~\ref{the:2.1}
and the inclusion $V_N\subseteq V_{N+1}$. If $V_{N+1}=V_N$, then
$p_iV_N\subseteq V_N$ for $i=1,2,3$. Since $1\in V_N$, repeated
left multiplication shows that every word belongs to $V_N$.
Conversely, if $\mathcal B=V_N$, then $V_{N+1}=V_N$.
In finite dimension the increasing sequence $(V_n)_n$ stabilizes.
\end{proof}

For each $x\in\mathcal M$, let $L_xy=xy$ for $y\in\mathcal M$.
The estimate $\|L_xy\|_2\leq\|x\|\|y\|_2$ gives a unique bounded
extension to $L^2(\mathcal M,\tau)$, also denoted by $L_x$.

\begin{theorem}\label{thm:layer-growth}
Let $H_k=V_k\ominus V_{k-1}$ and $h_k=\dim H_k$. Define
$A_{i,k}\colon H_k\to H_{k+1}$ by
$A_{i,k}\xi=(I-P_k)L_{p_i}\xi$. Then the following hold.
\begin{enumerate}
 \item For each $k\geq0$,
 \[
  H_{k+1}=\sum_{i=1}^3\operatorname{im}A_{i,k},
  \qquad A_{i,k+1}A_{i,k}=0,
 \]
 where the first sum is a linear sum.
 \item If $w$ is reduced of length $n+1$ and
 $r_w^{(n)}=(I-P_n)\widehat w$, then
 \[
  (I-P_{n+1})L_{p_i}\widehat w=A_{i,n+1}r_w^{(n)}.
 \]
 \item For $n\geq1$, one has $h_{n+1}\leq2h_n$ and
 $k_{n+1}\geq2k_n$. Consequently, for $n\geq0$,
 \[
  \frac{\operatorname{rank}S_{n+1}}{3\cdot 2^{n+1}}
  \leq\frac{\operatorname{rank}S_n}{3\cdot2^n},
  \qquad
  \frac{\dim\ker S_{n+1}}{3\cdot2^{n+1}}
  \geq\frac{\dim\ker S_n}{3\cdot2^n}.
 \]
\end{enumerate}
\end{theorem}

\begin{proof}
Since $V_{k+1}=V_k+\sum_{i=1}^3p_iV_k$ and
$p_iV_{k-1}\subseteq V_k$,
\[
 H_{k+1}=(I-P_k)V_{k+1}
 =\sum_{i=1}^3(I-P_k)L_{p_i}H_k
 =\sum_{i=1}^3\operatorname{im}A_{i,k}.
\]
For $\xi\in H_k$, the relation $p_i^2=p_i$ gives
\[
 \begin{aligned}
 A_{i,k+1}A_{i,k}\xi
 &=(I-P_{k+1})L_{p_i}(I-P_k)L_{p_i}\xi\\
 &=(I-P_{k+1})(p_i\xi-p_iP_kp_i\xi)=0,
 \end{aligned}
\]
since both terms inside the last parentheses belong to $V_{k+1}$.
For a reduced word of length $n+1$, write
$\widehat w=P_n\widehat w+r_w^{(n)}$. Then
$p_iP_n\widehat w\in V_{n+1}$, and hence
$(I-P_{n+1})L_{p_i}\widehat w=A_{i,n+1}r_w^{(n)}$.

For $k\geq1$, the inclusion
$\operatorname{im}A_{i,k-1}\subseteq\ker A_{i,k}$ gives
\[
 h_{k+1}\leq \sum_{i=1}^3\operatorname{rank}A_{i,k}
 \leq 3h_k-\sum_{i=1}^3\operatorname{rank}A_{i,k-1}
 \leq2h_k.
\]
The last inequality uses
$H_k=\sum_i\operatorname{im}A_{i,k-1}$. By Theorem~\ref{the:2.1},
\[
 k_{n+1}=3\cdot2^n-h_{n+1}
 \geq 2(3\cdot 2^{n-1}-h_n)=2k_n \qquad(n\geq1);
\]
Thus, applying Theorem~\ref{the:2.1} once more, we obtain for $n\geq0$,
 \[
  \frac{\operatorname{rank}S_{n+1}}{3\cdot 2^{n+1}}
  \leq\frac{\operatorname{rank}S_n}{3\cdot2^n},
  \qquad
  \frac{\dim\ker S_{n+1}}{3\cdot2^{n+1}}
  \geq\frac{\dim\ker S_n}{3\cdot2^n}.
 \]
\end{proof}

\begin{theorem}\label{thm:trace-decay}
With the notation of Theorem~\ref{the:2.1}, let $\operatorname{Tr}$
denote the unnormalized trace of a finite matrix. For every $n\geq0$,
\[
 \operatorname{Tr}S_{n+1}\leq\tfrac12\operatorname{Tr}S_n,\qquad
 \operatorname{Tr}S_0=\sum_{i=1}^3\tau(p_i)(1-\tau(p_i))\leq\tfrac34.
\]
Consequently,
\[
 0\leq\|S_n\|\leq\operatorname{Tr}S_n\leq\frac{3}{4\cdot2^n},
\]
and in particular $\|S_n\|\to0$ as $n\to\infty$.
\end{theorem}

\begin{proof}
Let $W_{n+1}$ denote the reduced formal words of length $n+1$.
For $w\in W_{n+1}$, write $r_w^{(n)}=(I-P_n)\widehat w\in V_{n+1}$.
If $i$ differs from the first letter of $w$, then $iw$ is reduced
of length $n+2$. Since $p_iP_n\widehat w\in p_iV_n\subseteq V_{n+1}$,
\[
 r_{iw}^{(n+1)}
 =(I-P_{n+1})L_{p_i}\widehat w
 =(I-P_{n+1})L_{p_i}P_{n+1}r_w^{(n)}.
\]
Set $b=(I-P_{n+1})L_{p_i}P_{n+1}$ and $a=P_{n+1}L_{p_i}P_{n+1}$.
Since $L_{p_i}$ is an orthogonal projection, $b^*b=a-a^2$.
The operator $a$ is a positive contraction, so functional
calculus gives $\|a-a^2\|\leq1/4$. Hence
$\|b\|=\|(I-P_{n+1})L_{p_i}P_{n+1}\|\leq1/2$ and
\[
 \|r_{iw}^{(n+1)}\|_2^2
 =\|(I-P_{n+1})L_{p_i}P_{n+1}r_w^{(n)}\|_2^2
 \leq\tfrac14\|r_w^{(n)}\|_2^2.
\]
Each $w$ admits two such prefixes, and every reduced word of length
$n+2$ arises once. It follows that
\begin{equation}\label{eq:layer-trace-step}
 \begin{aligned}
 \operatorname{Tr} S_{n+1}
 &=\sum_{w\in W_{n+1}}\sum_{\substack{1\leq i\leq 3\\i\ne\operatorname{first}(w)}}
   \|r_{iw}^{(n+1)}\|_2^2\\
 &\leq\sum_{w\in W_{n+1}}2\cdot\tfrac14\|r_w^{(n)}\|_2^2
   =\tfrac12\operatorname{Tr}S_n.
 \end{aligned}
\end{equation}

Since $V_0=\mathbb C1$, orthogonality gives
$0=\langle p_i-P_0p_i,1\rangle=\tau(p_i)-\tau(P_0p_i)$.
Thus $P_0p_i=\tau(p_i)1$ and
$p_i=\tau(p_i)1+(p_i-\tau(p_i)1)$. Since
\[
 \|p_i-\tau(p_i)1\|_2^2
 =\tau\!\left((p_i-\tau(p_i)1)^2\right)
 =\tau(p_i)-\tau(p_i)^2,
\]
we have
\begin{equation}\label{eq:layer-trace-start}
 \operatorname{Tr}S_0
 =\sum_{i=1}^3\bigl(\tau(p_i)-\tau(p_i)^2\bigr)
 =\frac34-\sum_{i=1}^3\left(\tau(p_i)-\frac12\right)^2
 \leq\frac34.
\end{equation}
By \eqref{eq:layer-trace-step}, \eqref{eq:layer-trace-start} and
positivity of $S_n$,
\[
 0\leq\|S_n\|\leq\operatorname{Tr}S_n
 \leq\frac{3}{4\cdot2^n},
 \qquad \|S_n\|\longrightarrow0.
\]
\end{proof}

\begin{remark}
The rank of $S_n$ equals $\dim(V_{n+1}\ominus V_n)$.
Theorem~\ref{thm:trace-decay} controls the squared norms of the
residuals. Finite-dimensionality is characterized by the
vanishing of some $S_N$, as in Corollary~\ref{cor:word-stabilization}.
In particular, $\mathcal B$ is finite-dimensional if
$\dim(V_{n+1}\ominus V_n)=1$ for all sufficiently large $n$.
The following theorem proves the stronger assertion that a single
one-dimensional layer suffices.
\end{remark}

\begin{theorem}\label{thm:one-layer}
Suppose $\dim H_d=1$ for some $d\geq1$. Fix $0\ne x\in H_d$ and let
$m(t)$ be the minimal polynomial of $x$, with $m=\deg m(t)$.
Then $V_{dm}=V_{dm-1}$, $\mathcal B=V_{dm-1}$, and $S_{dm-1}=0$.
Moreover, if $\mathcal B$ is infinite-dimensional, then
\[
 h_d\geq2\quad(d\geq1),\qquad
 \operatorname{rank}S_n\geq2\quad(n\geq0),\qquad
 \dim V_n\geq1+2n\quad(n\geq0).
\]
\end{theorem}

\begin{proof}
By hypothesis, $V_d=V_{d-1}+\mathbb Cx$. For arbitrary
$a_1,\ldots,a_m\in V_{d-1}+\mathbb Cx$, expand $a_1\cdots a_m$.
Every term except the multiple of $x^m$ contains a factor from
$V_{d-1}$ and hence has word length at most
$(m-1)d+(d-1)=dm-1$. Since $(V_d)^m=V_{dm}$, it follows that
\[
 V_{dm}\subseteq V_{dm-1}+\mathbb Cx^m.
\]
The minimal polynomial gives
\[
 x^m\in\operatorname{span}\{1,x,\ldots,x^{m-1}\}
 \subseteq V_{d(m-1)}\subseteq V_{dm-1}.
\]
Thus $V_{dm}=V_{dm-1}$, and Corollary~\ref{cor:word-stabilization}
gives $\mathcal B=V_{dm-1}$ and $S_{dm-1}=0$.

In infinite dimension the corollary excludes $h_d=0$, and the
argument above excludes $h_d=1$. The remaining assertions follow
from $\operatorname{rank}S_n=h_{n+1}$ and
$\dim V_n=1+\sum_{d=1}^n h_d$.
\end{proof}

\begin{theorem}\label{thm:uniform-minpoly}
For every $d\geq0$, there is an integer $D_d\geq1$ such that
$\deg m_x\leq D_d$ for every $x\in V_d$, where $m_x$ denotes the
minimal polynomial of $x$. Consequently,
$\#\sigma_{\mathcal M}(x)\leq D_d$ for every $x\in V_d$.
\end{theorem}

\begin{proof}
For each integer $m\geq1$, set
\[
 \begin{aligned}
 F_m&=\{x\in V_d:1,x,\ldots,x^m\text{ are linearly dependent}\}, \qquad
 \Gamma_m(x)&=\bigl(\tau((x^i)^*x^j)\bigr)_{i,j=0}^m.
 \end{aligned}
\]
The algebra $\mathcal B$ is algebraic, so $V_d=\bigcup_{m\geq1}F_m$.
The Gram matrix criterion gives $x\in F_m$ if and only if
$\det\Gamma_m(x)=0$. In particular, every $F_m$ is closed.
The Baire category theorem gives an $m_0$ such that $F_{m_0}$
has nonempty interior in $V_d$.

Choose a basis $e_1,\ldots,e_r$ of $V_d$ and for
$z=(z_1,\ldots,z_r)\in\mathbb C^r$ write
$x(z)=\sum_{k=1}^r z_ke_k$, with $z_k=x_k+iy_k$.
The function
\[
 g_m(x_1,y_1,\ldots,x_r,y_r)=\det\Gamma_m(x(z))
\]
is a polynomial in the real coordinates. We also use $g_m$ to
denote this function on $V_d$. Fix $b\in V_d$ and choose
$a\in\operatorname{int}(F_{m_0})$. The polynomial
$h_b(t)=g_{m_0}(a+t(b-a))$ in the real variable $t$ vanishes for
all sufficiently small $t$. Thus it vanishes identically. In
particular, $h_b(1)=g_{m_0}(b)=0$, so $b\in F_{m_0}$. Since $b$
was arbitrary, $F_{m_0}=V_d$ and we may take $D_d=m_0$.

To identify the spectrum, if $m_x(\lambda)\ne0$, then
 polynomial division gives
\[
 m_x(t)=(t-\lambda)r_\lambda(t)+m_x(\lambda),
 \qquad \deg r_\lambda=\deg m_x-1.
\]
Evaluating at \(x\) and using \(m_x(x)=0\), we obtain
\[
 (x-\lambda 1_{\mathcal M})r_\lambda(x)
 =r_\lambda(x)(x-\lambda 1_{\mathcal M})
 =-m_x(\lambda)1_{\mathcal M}.
\]
If \(m_x(\lambda)\ne0\), this identity gives the two-sided inverse
\[
 (x-\lambda 1_{\mathcal M})^{-1}
 =-\frac{r_\lambda(x)}{m_x(\lambda)}.
\]
If \(m_x(\lambda)=0\) and \(x-\lambda 1_{\mathcal M}\) were
invertible, the same identity would imply \(r_\lambda(x)=0\).
Since \(r_\lambda\) is a nonzero polynomial of smaller degree than
\(m_x\), this would contradict the minimality of \(m_x\). Hence
\[
 \sigma_{\mathcal M}(x)
 =\{\lambda\in\mathbb C:m_x(\lambda)=0\},
 \qquad
 \#\sigma_{\mathcal M}(x)\leq\deg m_x\leq D_d.
\]
\end{proof}

\begin{corollary}\label{cor:polarized-words}
Let $D=D_d$ be as in Theorem~\ref{thm:uniform-minpoly}.
For $x_1,\ldots,x_D\in V_d$, put
\[
 \operatorname{Sym}_D(x_1,\ldots,x_D)
 =\sum_{\pi\in\mathfrak S_D}x_{\pi(1)}\cdots x_{\pi(D)},
\]
where $\mathfrak S_D$ is the symmetric group on $\{1,\ldots,D\}$.
Then $\operatorname{Sym}_D(x_1,\ldots,x_D)\in V_{d(D-1)}$.
\end{corollary}

\begin{proof}
For $t=(t_1,\ldots,t_D)\in\mathbb C^D$ put
$x(t)=\sum_{j=1}^D t_jx_j\in V_d$. Let
$m_t(\lambda)=\lambda^s+\sum_{j=0}^{s-1}a_j\lambda^j$ be the
minimal polynomial of $x(t)$. By Theorem~\ref{thm:uniform-minpoly},
$s=\deg m_t\leq D$. Multiplying $m_t(x(t))=0$ by $x(t)^{D-s}$ yields
\[
 x(t)^D=-\sum_{j=0}^{s-1}a_jx(t)^{D-s+j}\in V_{d(D-1)},
\]
since $D-s+j\leq D-1$ for $0\leq j\leq s-1$.
Let $\rho\colon V_{dD}\to V_{dD}/V_{d(D-1)}$ be the quotient map.
The polynomial
\[
 F(t)=\rho(x(t)^D)
 =\rho\!\left(\sum_{i_1,\ldots,i_D=1}^D
     t_{i_1}\cdots t_{i_D}x_{i_1}\cdots x_{i_D}\right)
\]
vanishes identically. Taking the mixed partial derivative gives
\[
 0=\left.\frac{\partial^D F}{\partial t_1\cdots\partial t_D}
       \right|_{t=0}
 =\rho\bigl(\operatorname{Sym}_D(x_1,\ldots,x_D)\bigr).
\]
Thus $\operatorname{Sym}_D(x_1,\ldots,x_D)\in V_{d(D-1)}$.
\end{proof}

\begin{theorem}
Let $x$ be a nonzero element of $V_d$ with $d\geq1$, and let
$\{\lambda_1,\ldots,\lambda_r\}$ be the distinct nonzero spectral
values of $x^*x$. Then the right and left support projections
$r(x)$ and $l(x)$ belong to $V_{2dr}$.
\end{theorem}

\begin{proof}
Define the polynomial
\[
 f(t)=1-\prod_{i=1}^r\left(1-\frac{t}{\lambda_i}\right).
\]
It satisfies $f(0)=0$ and $f(\lambda_i)=1$ for every $i$.
Functional calculus gives $r(x)=f(x^*x)$. Since $V_d$ is
selfadjoint, $x^*x\in V_{2d}$, and $\deg f=r$ implies
$f(x^*x)\in V_{2dr}$. Likewise,
$\sigma(x^*x)\setminus\{0\}=\sigma(xx^*)\setminus\{0\}$ and
$xx^*\in V_{2d}$ give $l(x)=f(xx^*)\in V_{2dr}$.
\end{proof}

A first-return word begins and ends with $p_1$ and contains no other
occurrence of $p_1$. The intermediate letters of such a reduced word
alternate between $p_2$ and $p_3$. The following lemma embeds
these words in a fixed finite-dimensional linear space and
describes their role in the corner corresponding to $p_1$.

\begin{lemma}\label{lem:7}
Let $p=p_1$, $\mathcal C=\operatorname{alg}^*(1,p_2,p_3)$ and
$E=p\mathcal Cp$. Then $\mathcal C$ and $E$ are finite-dimensional.
The corner $\widetilde{\mathcal A}=p\mathcal Bp$ is generated by
$E$ as a $*$-algebra, and
\[
 \mathcal B=\mathcal C+\mathcal C\widetilde{\mathcal A}\mathcal C.
\]
If $\dim\mathcal C=c$ and $\dim\widetilde{\mathcal A}=a<\infty$,
then $\dim\mathcal B\leq c+c^2a$.
\end{lemma}

\begin{proof}
Theorem~\ref{thm: 2.4} gives $\dim\mathcal C<\infty$, and
$E$ is its image under compression by $p$. Put
$z=(p_2-p_3)^2=p_2-p_2p_3-p_3p_2+p_3$. Direct multiplication gives
\[
 \begin{gathered}
 zp_2=p_2z,\qquad zp_3=p_3z,\qquad
 p_3p_2=p_2+p_3-z-p_2p_3,\\
 p_2p_3p_2=(1-z)p_2,\qquad p_3p_2p_3=(1-z)p_3.
 \end{gathered}
\]
The first-return words fall into the four families
\[
 \begin{gathered}
 p_1(p_2p_3)^kp_1,\qquad p_1(p_3p_2)^kp_1 \qquad(k\geq1),\\
 p_1(p_2p_3)^kp_2p_1,\qquad p_1(p_3p_2)^kp_3p_1 \qquad(k\geq0).
 \end{gathered}
\]
Let $m(t)=t^s+\sum_{i=0}^{s-1}\alpha_i t^i$ be the minimal
polynomial of $p_2p_3$. For $k\geq0$ it gives
\[
 \begin{aligned}
 p_1(p_2p_3)^{k+s}p_1
 &=-\sum_{i=0}^{s-1}\alpha_i p_1(p_2p_3)^{k+i}p_1,\\
 p_1(p_2p_3)^{k+s}p_2p_1
 &=-\sum_{i=0}^{s-1}\alpha_i p_1(p_2p_3)^{k+i}p_2p_1.
 \end{aligned}
\]
The corresponding identities for the two $p_3p_2$ families
follow by using the polynomial with conjugated coefficients,
which annihilates $p_3p_2=(p_2p_3)^*$. Thus $E$ is spanned by
\[
 \begin{aligned}
 &\{p_1\}\,\cup\,
   \{p_1(p_2p_3)^kp_1,\ p_1(p_3p_2)^kp_1:1\leq k<s\}\\
 &\hspace{2em}\cup
   \{p_1(p_2p_3)^kp_2p_1,\ p_1(p_3p_2)^kp_3p_1:0\leq k<s\}.
 \end{aligned}
\]
In particular, $\dim E\leq 1+2(s-1)+2s=4s-1$.

Any word containing $p$ has the form $c_0pc_1p\cdots pc_k$ with
$k\geq1$ and $c_0,\ldots,c_k\in\mathcal C$. After compression
by $p$, such a word becomes a product of elements of $E$. In the original word, the portion from
the first to the last occurrence of $p$ belongs to the corner,
and the two outside portions belong to $\mathcal C$. Thus
$\mathcal B\subseteq\mathcal C+\mathcal C\widetilde{\mathcal A}\mathcal C$.
The reverse inclusion holds because $\mathcal C$ and
$\widetilde{\mathcal A}$ lie in $\mathcal B$. If the corner is
finite-dimensional, products of basis elements of the three
spaces give $\dim\mathcal B\leq c+c^2a$.
\end{proof}

\begin{theorem}\label{thm:corner-group}
Under the notation of Lemma~\ref{lem:7}, assume $p_1\ne0$.
There are an integer $N\geq2$ and finite unitaries
$u_1,\ldots,u_m\in\widetilde{\mathcal A}$, with
$m\leq\dim E-1$, such that
\[
 u_i^*u_i=u_iu_i^*=p_1,\qquad u_i^N=p_1,\qquad
 \widetilde{\mathcal A}=\operatorname{alg}^*(p_1,u_1,\ldots,u_m).
\]
If $\Gamma=\langle u_1,\ldots,u_m\rangle$ is the group they
generate with identity $p_1$, then
\[
 \widetilde{\mathcal A}=\operatorname{span}\Gamma
 \cong\mathbb C[\Gamma]/\ker\Pi,\qquad
 \Pi\left(\sum_{t\in\Gamma}c_t[t]\right)=\sum_{t\in\Gamma}c_tt.
\]
\end{theorem}

\begin{proof}
By Lemma~\ref{lem:7}, $E$ is finite-dimensional. Since $E$ is
closed under adjoints and contains $p_1$, it has a complex basis
$\{p_1,a_1,\ldots,a_m\}$ with $a_i=a_i^*$ and $m=\dim E-1$.
In the corner with unit $p_1$, write the spectral decompositions
\[
 a_i=\sum_{l=0}^{d_i}\lambda_{i,l}e_{i,l},\qquad
 \sum_{l=0}^{d_i}e_{i,l}=p_1.
\]
The eigenvalues $\lambda_{i,l}$ are distinct for each $i$.
Choose $N\geq\max\{2,d_1+1,\ldots,d_m+1\}$, set
$\xi=e^{2\pi i/N}$ and define
$u_i=\sum_{l=0}^{d_i}\xi^l e_{i,l}$. The orthogonal projections
$e_{i,l}$ give $u_i^*u_i=u_iu_i^*=p_1$ and $u_i^N=p_1$.

For each $i,l$, the Lagrange interpolation formula uses the polynomial
\[
 L_{i,l}(t)=\prod_{\substack{0\leq k\leq d_i\\k\ne l}}
 \frac{t-\lambda_{i,k}}{\lambda_{i,l}-\lambda_{i,k}}.
\]
Then
\[
 L_{i,l}(\lambda_{i,k})=
 \begin{cases}
  1,& k=l,\\
  0,& k\ne l
 \end{cases}
 \qquad(0\leq k\leq d_i).
\]
Thus $L_{i,l}(a_i)=e_{i,l}$, with constant terms evaluated using
the corner unit $p_1$, and
\[
 u_i=\sum_{l=0}^{d_i}\xi^l L_{i,l}(a_i)
 \in\operatorname{alg}^*(p_1,a_i).
\]
Since $\xi^0,\ldots,\xi^{d_i}$ are distinct, Lagrange interpolation at
these values gives a polynomial $\psi_i$ with
$\psi_i(\xi^l)=\lambda_{i,l}$. Hence
$a_i=\psi_i(u_i)\in\operatorname{alg}^*(p_1,u_i)$ and
\[
 \operatorname{alg}^*(p_1,a_1,\ldots,a_m)
 =\operatorname{alg}^*(p_1,u_1,\ldots,u_m)
 =\widetilde{\mathcal A}.
\]
The right side equals the linear span of $\Gamma$, so $\Pi$ is
a surjective $*$-homomorphism and the quotient formula follows.
\end{proof}

\begin{corollary}
Under the notation of Theorem~\ref{thm:corner-group}, if
$\ker\Pi=0$, then $\mathcal B$ is finite-dimensional.
\end{corollary}

\begin{proof}
If $\ker\Pi=0$, then $\widetilde{\mathcal A}\cong\mathbb C[\Gamma]$.
The corner $\widetilde{\mathcal A}$ is an $R^*$-algebra because
it is an algebraic $*$-subalgebra of $B(p_1\mathcal H)$ in a
faithful representation of $\mathcal M$ on $\mathcal H$.
By \cite[Example~3.9]{Mor}, the group algebra $\mathbb C[G]$
of a discrete group $G$ is an $R^*$-algebra exactly when $G$ is
locally finite. Hence $\Gamma$ is locally finite. It is also
finitely generated, so $\Gamma$ is finite and
$\widetilde{\mathcal A}$ is finite-dimensional.
Lemma~\ref{lem:7} gives the finite-dimensionality of $\mathcal B$.
\end{proof}

\section{A three-projection counterexample}\label{sec:counterexample}

In this section we use the Pauli-matrix form of the Jordan--Wigner construction
\cite{JordanWigner} to construct the counterexample we need.

\subsection{Clifford matrices and and the Vanishing of Cross Compressions}

\(\newline\)Let
\[
 p=\begin{pmatrix}0&1\\1&0\end{pmatrix},\qquad
 q=\begin{pmatrix}0&-i\\i&0\end{pmatrix},\qquad
 r=\begin{pmatrix}1&0\\0&-1\end{pmatrix}.
\]
Fix an integer $n\geq1$. Choose an even integer $g>n^2$,
write $g=2s$, and put $N=2^s$. On $(\mathbb C^2)^{\otimes s}$, define
\begin{equation}\label{eq: 4.1}
 \begin{split}
 \gamma_{2t-1}&=r^{\otimes(t-1)}\otimes p\otimes I^{\otimes(s-t)},\\
 \gamma_{2t}&=r^{\otimes(t-1)}\otimes q\otimes I^{\otimes(s-t)}
 \qquad(1\leq t\leq s).
 \end{split}
\end{equation}
The Pauli matrices are selfadjoint involutions and pairwise anticommute.
Generators at distinct positions anticommute because their tensor factors
anticommute at the earlier position and commute at the other positions.
It follows that
\begin{equation}\label{eq: 4.2}
 \gamma_i^*=\gamma_i,\qquad
 \gamma_i^2=I_N,\qquad
 \gamma_i\gamma_j=-\gamma_j\gamma_i\quad(i\ne j).
\end{equation}
Since $pq=ir$, multiplying the two generators at position $t$ gives
$iI^{\otimes(t-1)}\otimes r\otimes I^{\otimes(s-t)}$. Consequently,
\begin{equation}\label{eq: 4.3}
 \Gamma:=r^{\otimes s}=(-i)^s\gamma_1\cdots\gamma_g.
\end{equation}
The matrix $\Gamma$ is a real diagonal selfadjoint involution with
$\operatorname{Tr}\Gamma=(\operatorname{Tr}r)^s=0$. Thus its diagonal
entries $+1$ and $-1$ each occur $N/2$ times.

\begin{lemma}\label{lem:clifford}
For $I=\{i_1<\cdots<i_r\}\subseteq\{1,\ldots,g\}$, put
$\gamma_I=\gamma_{i_1}\cdots\gamma_{i_r}$ and
$\gamma_{\varnothing}=I_N$. The family $\{\gamma_I\}$ is a
trace-orthogonal basis of $M_N(\mathbb C)$.
If $W$ is a complex subspace of
$\operatorname{span}_{\mathbb C}\{\gamma_1,\ldots,\gamma_g\}$ and
$h=\dim_{\mathbb C}W<g$, then the unital algebra $\mathcal A(W)$
generated by $W$ satisfies $\operatorname{Tr}(\Gamma a)=0$ for
every $a\in\mathcal A(W)$.
\end{lemma}

\begin{proof}
The identities \eqref{eq: 4.2} reduce any word in the $\gamma_i$
to $\pm\gamma_I$ by interchanging distinct generators and cancelling
repeated ones. We first show that $\operatorname{Tr}(\gamma_I)=0$
for $I\ne\varnothing$. Set $m=|I|>0$. If $m$ is even, choose
$j\in I$. Then
$\gamma_j\gamma_I\gamma_j=(-1)^{m-1}\gamma_I=-\gamma_I$,
so by the cyclicity of the trace,
$\operatorname{Tr}(\gamma_I)=\operatorname{Tr}(\gamma_j\gamma_I\gamma_j)
=\operatorname{Tr}(-\gamma_I)$.
If $m$ is odd, the evenness of $g$ gives some $j\notin I$, and
$\gamma_j\gamma_I\gamma_j=(-1)^m\gamma_I=-\gamma_I$.
A similar calculation then yields $\operatorname{Tr}(\gamma_I)=0$.

If $I\ne J$, then $\gamma_I^*\gamma_J=\pm\gamma_{I\triangle J}$,
where $I\triangle J\ne\varnothing$. Thus
$\operatorname{Tr}(\gamma_I^*\gamma_J)=
\operatorname{Tr}(\pm\gamma_{I\triangle J})=0$.
For $I=J$, unitarity gives
$\operatorname{Tr}(\gamma_I^*\gamma_I)=\operatorname{Tr}(I_N)=N$.
Consequently, $N^{-1}\operatorname{Tr}(\gamma_I^*\gamma_J)=\delta_{I,J}$.
There are $2^g=N^2$ subsets of $\{1,\ldots,g\}$, so the
corresponding nonzero orthogonal matrices form a basis of
$M_N(\mathbb C)$, which has dimension $N^2$.

For $x=\sum_i \lambda_i\gamma_i$ and $y=\sum_i \mu_i\gamma_i$, define
\begin{equation}\label{eq: 4.4}
\beta(x,y)=\sum_i\lambda_i\mu_i.
\end{equation}
The identities \eqref{eq: 4.2} imply $xy+yx=2\beta(x,y)I_N$.
Choose a basis $w_1,\ldots,w_h$ of $W$. Then
\begin{equation}\label{eq: 4.5}
 w_jw_i=-w_iw_j+2\beta(w_i,w_j)I_N,\qquad
 w_i^2=\beta(w_i,w_i)I_N.
\end{equation}
We show that $\mathcal A(W)$ is spanned by $I_N$ and the products
$w_{j_1}\cdots w_{j_r}$ with $1\leq j_1<\cdots<j_r\leq h$.
Argue by induction on word length and then on inversion number for
words of a fixed length. Interchanging an adjacent inversion using
\eqref{eq: 4.5} decreases the inversion number and produces an
additional word of length smaller by two. Once the indices are
nondecreasing, apply the square relation to eliminate adjacent
repetitions. The resulting words have strictly increasing indices.
The argument also applies when $\beta$ is degenerate on $W$.

Let $F_m=\operatorname{span}_{\mathbb C}\{\gamma_I:|I|\leq m\}$.
Each $w_j$ is a linear combination of the $\gamma_i$, so
$w_{j_1}\cdots w_{j_r}$ expands into words of length $r$ in the
$\gamma_i$. The identities \eqref{eq: 4.2} never increase length,
hence $w_{j_1}\cdots w_{j_r}\in F_r$ and
$\mathcal A(W)\subseteq F_h$.

Write $a=\sum_{|I|\leq h}c_I\gamma_I$ for $a\in\mathcal A(W)$
and let $G=\{1,\ldots,g\}$. By \eqref{eq: 4.3},
$\Gamma\gamma_I=\eta_I\gamma_{G\setminus I}$ with $|\eta_I|=1$.
Since $h<g$, every complement $G\setminus I$ is nonempty. Hence
\[
 \operatorname{Tr}(\Gamma a)=
 \sum_{|I|\leq h}c_I\eta_I\operatorname{Tr}(\gamma_{G\setminus I})=0.
\]

\end{proof}

Let $V_g=\{0,a_1,b_1,\ldots,a_g,b_g\}$ and
$\mathbb C^{V_g}=\{\sum_{v\in V_g}\xi_ve_v:\xi_v\in\mathbb C\}$,
with orthonormal basis $(e_v)_{v\in V_g}$. Denote the corresponding
matrix units by $E_{s,t}$ for $s,t\in V_g$. On
$\mathbb C^{V_g}\otimes\mathbb C^N$ define
\begin{equation}\label{eq:4.6}
 K^{\mathrm{bouq}}=\sum_{i=1}^g\left[
 (E_{0a_i}+E_{a_i0}+E_{a_ib_i}+E_{b_ia_i})\otimes I_N
 +(E_{0b_i}+E_{b_i0})\otimes\gamma_i\right].
\end{equation}
For each $i$, the summand
\[
 (E_{0a_i}+E_{a_i0}+E_{a_ib_i}+E_{b_ia_i})\otimes I_N
 +(E_{0b_i}+E_{b_i0})\otimes\gamma_i
\]
vanishes on the orthogonal complement of
$(e_0\otimes\mathbb C^N)\oplus(e_{a_i}\otimes\mathbb C^N)
\oplus(e_{b_i}\otimes\mathbb C^N)$. Write
$K^{\mathrm{bouq}}=\sum_{i=1}^g T_i$, where $T_i$ is the extension
by zero of the matrix
\[
 \begin{pmatrix}
 0&I_N&\gamma_i\\
 I_N&0&I_N\\
 \gamma_i&I_N&0
 \end{pmatrix}
\]
on that subspace. For $x,y,z\in\mathbb C^N$, its action there gives
\[
 \begin{aligned}
 \|T_i(x,y,z)^T\|^2
 &=\|y+\gamma_i z\|^2+\|z+x\|^2+\|\gamma_ix+y\|^2\\
 &\leq 4\bigl(\|x\|^2+\|y\|^2+\|z\|^2\bigr).
 \end{aligned}
\]
Hence $\|T_i\|\leq2$ and $\|K^{\mathrm{bouq}}\|\leq2g$.

Set
\[
 \mathcal K=\bigoplus_{c=1}^3
  (\mathbb C^{V_g}\otimes\mathbb C^N\otimes\mathbb C^N),
 \qquad R=\dim\mathcal K=3(2g+1)N^2.
\]
Define $\rho:V_g\to\{0,1,\ldots,2g\}$ by
\[
 \rho(0)=0,\qquad \rho(a_i)=2i-1,\qquad \rho(b_i)=2i.
\]
Let $L=3(2g+1)$ and put
\[
 \nu_{c,v}=\frac14+
 \frac{(c-1)(2g+1)+\rho(v)+1}{2(L+1)}.
\]
The numbers $\nu_{c,v}$ are pairwise distinct and lie in $(1/4,3/4)$.
Define
\[
 \begin{aligned}
 H&=\bigoplus_{c=1}^3
       \left(\sum_{v\in V_g}\nu_{c,v}E_{vv}\otimes I_N\otimes I_N\right),\\
 K&=\frac12I_R+\frac1{8g}\bigoplus_{c=1}^3
        (K^{\mathrm{bouq}}\otimes I_N).
 \end{aligned}
\]
The bounds on the vertex markers and on $\|K^{\mathrm{bouq}}\|$
give $\frac14I_R<H<\frac34I_R$ and $\frac14I_R\leq K\leq\frac34I_R$.

\begin{lemma}\label{lem:4.2}
Let $d=(2g+1)N$, and set $B_0:=\bigoplus_{c=1}^3
          (M_d(\mathbb{C})\otimes I_N).$ Then
${\mathcal{A}}^*(1,H,K)=B_0.$
\end{lemma}

\begin{proof}
Let $D=\operatorname{alg}^{*}(1,H,K)$. Since $1,H,K\in B_0$,
it follows that $D\subseteq B_0$. Hence it suffices to show that $B_0\subseteq D$.

For each $u\in V_g$, the Lagrange interpolation formula gives
\[
 Q_{1,u}(t)=
 \prod_{\substack{(d,w)\in\{1,2,3\}\times V_g\\(d,w)\ne(1,u)}}
 \frac{t-\nu_{d,w}}{\nu_{1,u}-\nu_{d,w}}.
\]
Its value at $H$ is the spectral projection corresponding to $\nu_{1,u}$.
For brevity, write $Q_{1,u}=Q_{1,u}(H)$. Set
\[
 P_u=E_{uu}\otimes I_N\otimes I_N,\qquad
 \widetilde K=\frac12I_{(2g+1)N^2}
             +\frac1{8g}(K^{\mathrm{bouq}}\otimes I_N).
\]
Then
\[
 Q_{1,u}
 =
 \begin{pmatrix}
 P_u&0&0\\
 0&0&0\\
 0&0&0
 \end{pmatrix},
 K
 =
 \begin{pmatrix}
 \widetilde K&0&0\\
 0&\widetilde K&0\\
 0&0&\widetilde K
 \end{pmatrix}.
\]
Block multiplication gives
\[
 Q_{1,u}KQ_{1,v}
 =
 \begin{pmatrix}
 P_u\widetilde KP_v&0&0\\
 0&0&0\\
 0&0&0
 \end{pmatrix}.
\]

Using
$E_{uu}E_{\alpha\beta}E_{vv}=\delta_{u,\alpha}\delta_{\beta,v}E_{uv}
$
and the definition of $K^{\mathrm{bouq}}$, we obtain
\[
 Q_{1,u}KQ_{1,v}
 =
 \begin{cases}
 \displaystyle \frac12Q_{1,u},
 &u=v,\\[5pt]
 \displaystyle
 \frac1{8g}(E_{uv}\otimes I_N\otimes I_N)\oplus0\oplus0,
 &\{u,v\}=\{0,a_i\}\ \text{or}\ \{a_i,b_i\},\\[5pt]
 \displaystyle
 \frac1{8g}(E_{uv}\otimes\gamma_i\otimes I_N)\oplus0\oplus0,
 &\{u,v\}=\{0,b_i\},\\[5pt]
 0,&\text{otherwise},
 \end{cases}
\]
for $1\le i\le g$.

For $u,v\in V_g$, write
$F_{uv}:=
 (E_{uv}\otimes I_N\otimes I_N)\oplus0\oplus0.$
The preceding calculation shows that both
$F_{0a_i}=8gQ_{1,0}KQ_{1,a_i}$ and $
 F_{a_ib_i}=8gQ_{1,a_i}KQ_{1,b_i}$ are in $D$. Since
 $F_{0b_i}
 =F_{0a_i}F_{a_ib_i}
 =(E_{0a_i}E_{a_ib_i}\otimes I_N\otimes I_N)
      \oplus0\oplus0=(E_{0b_i}\otimes I_N\otimes I_N)\oplus0\oplus0
 \in D,$
it follows that $F_{u0},F_{0u}\in D$ for every $u\in V_g$, including
$F_{00}=Q_{1,0}$. Consequently,
\[
 F_{uv}=F_{u0}F_{0v}\in D
(u,v\in V_g).
\]

The preceding calculation also gives $8gQ_{1,0}KQ_{1,b_i}=(E_{0b_i}\otimes\gamma_i\otimes I_N)\oplus0\oplus0.$
Multiplying by $F_{b_i0}$, we obtain
\[
 \begin{aligned}
 C_i
 &:=(8gQ_{1,0}KQ_{1,b_i})F_{b_i0}\\
 &=(E_{0b_i}E_{b_i0}\otimes\gamma_i\otimes I_N)
      \oplus0\oplus0\\
 &=(E_{00}\otimes\gamma_i\otimes I_N)\oplus0\oplus0
 \in D.
 \end{aligned}
\]
Moreover,
$C_{i_1}\cdots C_{i_r}
 =(E_{00}\otimes\gamma_{i_1}\cdots\gamma_{i_r}\otimes I_N)
 \oplus0\oplus0.$
By Lemma~\ref{lem:clifford}, these products, together with $F_{00}$
for the empty product, span
\[
 (E_{00}\otimes M_N(\mathbb C)\otimes I_N)\oplus0\oplus0.
\]
Hence, for every $A\in M_N(\mathbb C)$,
$(E_{00}\otimes A\otimes I_N)\oplus0\oplus0\in D.$

Since
\[
 \begin{aligned}
F_{u0}
 \bigl[(E_{00}\otimes A\otimes I_N)\oplus0\oplus0\bigr]
 F_{0v}
 &=(E_{u0}E_{00}E_{0v}\otimes A\otimes I_N)\oplus0\oplus0\\
  &=(E_{uv}\otimes A\otimes I_N)\oplus0\oplus0
 \in D.
 \end{aligned}
\]
it follows that
$(M_{2g+1}(\mathbb C)\otimes M_N(\mathbb C)\otimes I_N)\oplus0\oplus0
\subseteq D$. The same calculation in the other two summands gives
$B_0\subseteq D$.

\end{proof}

Write
\[
 \mathcal K=\bigoplus_{c=1}^3\mathcal K_c,\qquad
 \mathcal K_c=\mathbb C^{V_g}\otimes\mathbb C^N\otimes\mathbb C^N.
\]
Let $e_0$ be the basis vector of $\mathbb C^{V_g}$ at the vertex $0$,
and let $e_1,\ldots,e_N$ be the standard orthonormal basis of
$\mathbb C^N$. Set
\begin{equation}\label{eq: 4.7}
 \Omega=\frac1{\sqrt N}\sum_{a=1}^Ne_a\otimes e_a,\qquad
 u=e_0\otimes\Omega,\qquad
 w=e_0\otimes(\Gamma\otimes I_N)\Omega,
\end{equation}
For each $A\in M_N(\mathbb C)$,
\[
 \begin{aligned}
 \Omega^*(A\otimes I_N)\Omega
 &=\frac1N\sum_{a,b=1}^N(e_a^*Ae_b)(e_a^*e_b)\\
 &=\frac1N\sum_{a=1}^Ne_a^*Ae_a
 =\frac1N\operatorname{Tr}A.
 \end{aligned}
\]

Since $\Gamma^*=\Gamma$, $\Gamma^2=I_N$, and
$\operatorname{Tr}\Gamma=0$, it follows that
$u^*u=w^*w=1,~u^*w=0.$
With respect to the orthogonal decomposition $\mathcal K$,
set
\[
 \begin{aligned}
 u_1&=\begin{pmatrix}u\\0\\0\end{pmatrix},
 &u_2&=\begin{pmatrix}0\\u\\0\end{pmatrix},
 &u_3&=\begin{pmatrix}0\\0\\u\end{pmatrix},\\
 w_1&=\begin{pmatrix}w\\0\\0\end{pmatrix},
 &w_2&=\begin{pmatrix}0\\w\\0\end{pmatrix},
 &w_3&=\begin{pmatrix}0\\0\\w\end{pmatrix}.
 \end{aligned}
\]
Thus
$u_c^*u_d=w_c^*w_d=\delta_{cd}, ~u_c^*w_d=0 ~(1\le c,d\le3).$

Let $T$ be an operator on $\mathcal K$ of the form
\[
 T=[\widehat T_{cd}\otimes I_N]_{c,d=1}^3.
\]
In full block form this is
\[
 T=
 \begin{pmatrix}
 \widehat T_{11}\otimes I_N&
 \widehat T_{12}\otimes I_N&
 \widehat T_{13}\otimes I_N\\
 \widehat T_{21}\otimes I_N&
 \widehat T_{22}\otimes I_N&
 \widehat T_{23}\otimes I_N\\
 \widehat T_{31}\otimes I_N&
 \widehat T_{32}\otimes I_N&
 \widehat T_{33}\otimes I_N
 \end{pmatrix},
\]
where $\widehat T_{cd}\otimes I_N:{\mathcal{K}}_d\to {\mathcal{K}}_c$ for $c,d\in \{1,2,3\}$. For a fixed $c$, decompose
$\widehat T_{cc}$ according to the basis
$\{e_0,e_{a_1},e_{b_1},\ldots,e_{a_g},e_{b_g}\}$ of $\mathbb C^{V_g}$:
\[
 \widehat T_{cc}
 =\bigl[T_{uv}^{(c)}\bigr]_{u,v\in V_g},\qquad  T_{uv}^{(c)}\in M_N(\mathbb C).
\]

For $c=1$, block multiplication gives
\[
 u_1^*Tu_1=u^*(\widehat T_{11}\otimes I_N)u.
\]
Since $u=e_0\otimes\Omega$, its column is
$(\Omega,0,\ldots,0)^T$. Consequently,
\[
 \begin{aligned}
 u_1^*Tu_1
 &=
 \begin{pmatrix}\Omega^*&0&\cdots\end{pmatrix}
 \begin{pmatrix}
 T_{00}^{(1)}\otimes I_N&
 T_{0a_1}^{(1)}\otimes I_N&\cdots\\
 T_{a_10}^{(1)}\otimes I_N&
 T_{a_1a_1}^{(1)}\otimes I_N&\cdots\\
 \vdots&\vdots&\ddots
 \end{pmatrix}
 \begin{pmatrix}\Omega\\0\\\vdots\end{pmatrix}\\
 &=\Omega^*(T_{00}^{(1)}\otimes I_N)\Omega
 =\frac1N\operatorname{Tr}T_{00}^{(1)}.
 \end{aligned}
\]

The same computation applies to every $c$. Replacing either copy of $\Omega$ by $(\Gamma\otimes I_N)\Omega$ yields
\begin{equation}\label{eq:4.8}
 \begin{aligned}
 u_c^*Tu_c
 &=\frac1N\operatorname{Tr}T_{00}^{(c)},\\
 u_c^*Tw_c
 &=\frac1N\operatorname{Tr}(T_{00}^{(c)}\Gamma),\\
 w_c^*Tu_c
 &=\frac1N\operatorname{Tr}(\Gamma T_{00}^{(c)}),\\
 w_c^*Tw_c
 &=\frac1N\operatorname{Tr}(\Gamma T_{00}^{(c)}\Gamma)
 =\frac1N\operatorname{Tr}T_{00}^{(c)}.
 \end{aligned}
\end{equation}

These identities apply entry by entry after adding an auxiliary factor.
Suppose $T$ acts on $\mathbb C^k\otimes\mathcal K$ and that each
auxiliary matrix entry has the block form above. Let
\[
 T_{00}^{(c)}\in M_k(\mathbb C)\otimes M_N(\mathbb C)
\]
be the matrix of the corresponding root--root blocks, and put
$\tau_N=N^{-1}\operatorname{Tr}$. Then
\begin{equation}\label{eq:4.9}
 \begin{aligned}
 (I_k\otimes u_c)^*T(I_k\otimes u_c)
 &=(\operatorname{id}_{M_k}\otimes\tau_N)(T_{00}^{(c)}),\\
 (I_k\otimes u_c)^*T(I_k\otimes w_c)
 &=(\operatorname{id}_{M_k}\otimes\tau_N)
       \bigl(T_{00}^{(c)}(I_k\otimes\Gamma)\bigr),\\
 (I_k\otimes w_c)^*T(I_k\otimes u_c)
 &=(\operatorname{id}_{M_k}\otimes\tau_N)
       \bigl((I_k\otimes\Gamma)T_{00}^{(c)}\bigr),\\
 (I_k\otimes w_c)^*T(I_k\otimes w_c)
 &=(\operatorname{id}_{M_k}\otimes\tau_N)(T_{00}^{(c)}).
 \end{aligned}
\end{equation}

\begin{lemma}\label{lem: 4.3}
Let $1\leq k^2\leq n^2 <g$ and $Z_0,X_0,Y_0\in M_k(\mathbb{C})$. If
\[
 \Delta=Z_0\otimes I_R-X_0\otimes H-Y_0\otimes K
\]
is invertible, then $(I_k\otimes u_c)^*\Delta^{-1}(I_k\otimes w_c)=(I_k\otimes w_c)^*\Delta^{-1}(I_k\otimes u_c)=0$
for $c=1,2,3$. The coefficients need not commute or be Hermitian.
\end{lemma}

\begin{proof}
Since $H$ and $K$ preserve each summand and act trivially on
the last tensor factor, we may write
\[
 \Delta =\oplus_{c=1}^3(\widehat\Delta_c\otimes I_N),
\]
where $\widehat\Delta_c$ acts on
$\mathcal L_c=\mathbb C^k\otimes\mathbb C^{V_g}\otimes\mathbb C^N$.
In detail,
\[
 \begin{aligned}
 \widehat\Delta_c
 &=Z_0\otimes I_{(2g+1)N}
   -X_0\otimes\left(\sum_{v\in V_g}\nu_{c,v}E_{vv}\otimes I_N\right)
 -Y_0\otimes\left(\frac12 I_{(2g+1)N}
              +\frac1{8g}K^{\mathrm{bouq}}\right).
 \end{aligned}
\]
The invertibility of $\Delta$ implies that each $\widehat\Delta_c$
is invertible. Fix $c$ and use the orthogonal decomposition
\[
 \mathbb C^{V_g}=\operatorname{span}\{e_0\}\oplus
    \bigoplus_{i=1}^g\operatorname{span}\{e_{a_i},e_{b_i}\}.
\]
This gives
\[
 \mathcal L_c=\mathcal E_0\oplus(\oplus_{i=1}^g\mathcal E_i),
\]
where
\[
 \begin{aligned}
 \mathcal E_0=\mathbb C^k\otimes\operatorname{span}\{e_0\}
        \otimes\mathbb C^N,\qquad
 \mathcal E_i=\mathbb C^k\otimes\operatorname{span}\{e_{a_i},e_{b_i}\}
        \otimes\mathbb C^N.
 \end{aligned}
\]
Thus $\dim\mathcal E_0=kN$ and $\dim\mathcal E_i=2kN$.

Put $A=Y_0/(8g)$ and $D_v=Z_0-\nu_{c,v}X_0-\frac12Y_0$.
For every $v\in V_g$, we define
\[
J_v:\mathbb C^k\otimes\mathbb C^N
\longrightarrow
\mathbb C^k\otimes\mathbb C^{V_g}\otimes\mathbb C^N,
\qquad
J_v(\eta\otimes\xi)=\eta\otimes e_v\otimes\xi.
\]
For $u,v\in V_g$, put $\Delta_{uv}:=J_u^*\widehat\Delta_cJ_v$. It follows that
\[
 \Delta_{vv}
 =Z_0\otimes I_N-\nu_{c,v}X_0\otimes I_N
       -\frac12Y_0\otimes I_N
 =D_v\otimes I_N
\]
and
\[
 \begin{aligned}
 \Delta_{0a_i}=\Delta_{a_i0}&=-A\otimes I_N,\\
 \Delta_{a_ib_i}=\Delta_{b_ia_i}&=-A\otimes I_N,\\
 \Delta_{0b_i}=\Delta_{b_i0}&=-A\otimes\gamma_i.
 \end{aligned}
\]
All other off-diagonal blocks vanish. Consequently, the principal block indexed by $0,a_i,b_i$ is
\[
 \begin{pmatrix}
 D_0\otimes I_N&-A\otimes I_N&-A\otimes\gamma_i\\
 -A\otimes I_N&D_{a_i}\otimes I_N&-A\otimes I_N\\
 -A\otimes\gamma_i&-A\otimes I_N&D_{b_i}\otimes I_N
 \end{pmatrix}.
\]
Define
\[
 \mathcal D_i=
 \begin{pmatrix}
 D_{a_i}&-A\\
 -A&D_{b_i}
 \end{pmatrix},
\quad \mathcal R_i=
 \begin{pmatrix}-A\otimes I_N&-A\otimes\gamma_i\end{pmatrix},
~
 \mathcal C_i=
 \begin{pmatrix}-A\otimes I_N\\-A\otimes\gamma_i\end{pmatrix}.
\]
Using $\mathcal L_c=\mathcal E_0\oplus\left(\oplus_{i=1}^g\mathcal E_i\right)$, we obtain that
\[
 \widehat\Delta_c=
 \begin{pmatrix}
 D_0\otimes I_N&\mathcal R_1&\cdots&\mathcal R_g\\
 \mathcal C_1&\mathcal D_1\otimes I_N&\cdots&0\\
 \vdots&\vdots&\ddots&\vdots\\
 \mathcal C_g&0&\cdots&\mathcal D_g\otimes I_N
 \end{pmatrix}.
\]
Since the auxiliary coefficients are not assumed selfadjoint,
$\mathcal C_i$ need not equal $\mathcal R_i^*$.

First we assume that all branch matrices $\mathcal D_i$ are
invertible. Write
\[
 \mathcal D_i^{-1}
 =
 \begin{pmatrix}
 E_i&F_i\\
 G_i&L_i
 \end{pmatrix},
\]
for some $E_i,F_i,G_i,L_i\in M_k(\mathbb C)$.
Direct multiplication yields
\[
 \begin{aligned}
 \mathcal R_i(\mathcal D_i^{-1}\otimes I_N)\mathcal C_i
 &=AE_iA\otimes I_N+AF_iA\otimes\gamma_i
 +AG_iA\otimes\gamma_i+AL_iA\otimes\gamma_i^2\\
 &=A(E_i+L_i)A\otimes I_N+A(F_i+G_i)A\otimes\gamma_i.
 \end{aligned}
\]
The Schur complement $S$ of the branch diagonal block in
$\widehat\Delta_c$ is invertible and has the form
\[
 S=C_0\otimes I_N+\sum_{i=1}^gZ_i\otimes\gamma_i,
\]
where
\[
 C_0=D_0-\sum_{i=1}^g A(E_i+L_i)A,\qquad Z_i=-A(F_i+G_i)A.
\]
Write
\[
\widehat\Delta_c^{-1}
=
\begin{pmatrix}
{B_{00}}&B_{01}&\cdots&B_{0g}\\
B_{10}&B_{11}&\cdots&B_{1g}\\
\vdots&\vdots&\ddots&\vdots\\
B_{g0}&B_{g1}&\cdots&B_{gg}
\end{pmatrix}.
\]
Then the Schur complement formula implies that $B_{00}$, the upper-left $kN\times kN$ block of
$\widehat\Delta_c^{-1}$, is $S^{-1}$.

For $1\le a,b\le k$, let
\[
 \zeta_{ab}=\sum_{i=1}^g(Z_i)_{ab}\gamma_i,\qquad
 W=\operatorname{span}_{\mathbb C}
       \{\zeta_{ab}:1\le a,b\le k\},
\]
where $(Z_i)_{ab}$ denotes the element in the $a$-th row and $b$-th column of the matrix $Z_i$.
Then
$\dim_{\mathbb C}W\le k^2\le n^2<g.$

Let $\mathcal A(W)$ be the unital algebra generated by $W$.
Since $S_{ab}=(C_0)_{ab}I_N+\zeta_{ab}\in\mathcal A(W)$, it follows that
$ S\in M_k(\mathcal A(W)).$

Let $m=kN$ and
$P(t)=t^m+c_{m-1}t^{m-1}+\cdots+c_1t+c_0$ be the
characteristic polynomial of $S$. The invertibility of $S$ gives
$c_0\ne0$. By the Cayley--Hamilton theorem,
\[
 P(S)=S^m+c_{m-1}S^{m-1}+\cdots+c_1S+c_0I_m=0.
\]
It follows that
\[
 S^{-1}=-\frac1{c_0}
     (S^{m-1}+c_{m-1}S^{m-2}+\cdots+c_1I_m)
 \in M_k(\mathcal A(W)).
\]
By Lemma~\ref{lem:clifford},
\[
 \operatorname{Tr}\bigl(\Gamma(S^{-1})_{ab}\bigr)=0
 (1\le a,b\le k).
\]

Let $\tau_N=N^{-1}\operatorname{Tr}$. Applying \eqref{eq:4.9}
to $\Delta^{-1}$ and using cyclicity of the trace gives
\[
 \begin{aligned}
 (I_k\otimes u_c)^*\Delta^{-1}(I_k\otimes w_c)
 &=(\operatorname{id}_{M_k}\otimes\tau_N)
       (S^{-1}(I_k\otimes\Gamma))
 =0,\\
 (I_k\otimes w_c)^*\Delta^{-1}(I_k\otimes u_c)
 &=
 (\operatorname{id}_{M_k}\otimes\tau_N)
       ((I_k\otimes\Gamma)S^{-1})
 =0.
 \end{aligned}
\]
Since $c$ is arbitrary, these identities hold for all three
Hilbert space summands whenever the branch matrices are invertible.

Finally, we remove the additional branch invertibility assumption.
Fix $Z_0,X_0,Y_0$ for which $\Delta$ is invertible. For every $t\in\mathbb C$, define
\[
 Z_0(t)=Z_0+tI_k,
 ~\Delta(t)
 =Z_0(t)\otimes I_R-X_0\otimes H-Y_0\otimes K
 =\Delta+tI_{kR}.
\]
For each fixed $c$, the corresponding matrix on $\mathcal L_c$ is
\[
 \widehat\Delta_c(t)=
 \begin{pmatrix}
 (D_0+tI_k)\otimes I_N
     &\mathcal R_1&\cdots&\mathcal R_g\\
 \mathcal C_1
     &(\mathcal D_1+tI_{2k})\otimes I_N&\cdots&0\\
 \vdots&\vdots&\ddots&\vdots\\
 \mathcal C_g
     &0&\cdots&(\mathcal D_g+tI_{2k})\otimes I_N
 \end{pmatrix}.
\]
Thus the off-diagonal blocks are unchanged, while each branch
matrix becomes
\[
 \mathcal D_i(t)
 =
 \begin{pmatrix}
 D_{a_i}+tI_k&-A\\
 -A&D_{b_i}+tI_k
 \end{pmatrix}
 =\mathcal D_i+tI_{2k}.
\]

Put $p(t):=\det(\Delta+tI_{kR})$ and $q_{c,i}(t):=\det(\mathcal D_i^{(c)}+tI_{2k})$, where the superscript $(c)$ labels the
Hilbert space summand.
Both determinant polynomials are monic, of degrees $kR$ and
$2k$, respectively.
\[
 \mathcal F
 =
 \{t:p(t)=0\}
 \cup
 \bigcup_{c=1}^3\bigcup_{i=1}^g\{t:q_{c,i}(t)=0\}.
\]
Then, for every $t\notin\mathcal F$, both $\Delta(t)$ and all branch
matrices are invertible. Therefore, the preceding argument gives
\[
 \begin{aligned}
 (I_k\otimes u_c)^*\Delta(t)^{-1}(I_k\otimes w_c)=0,\qquad
 (I_k\otimes w_c)^*\Delta(t)^{-1}(I_k\otimes u_c)=0.
 \end{aligned}
\]

Define
\[
 \begin{aligned}
 \mathcal P_c(t)
 &= (I_k\otimes u_c)^*\operatorname{adj}(\Delta(t))(I_k\otimes w_c),\\
 \mathcal Q_c(t)
 &= (I_k\otimes w_c)^*\operatorname{adj}(\Delta(t))(I_k\otimes u_c).
 \end{aligned}
\]
Then both are in $M_k(\mathbb C[t])$.

For $t\notin\mathcal F$, we get that $p(t)=\det\Delta(t)\ne0$ and
$\Delta(t)^{-1}=\frac{\operatorname{adj}(\Delta(t))}{p(t)}$.
Hence
\[
 \begin{aligned}
 \mathcal P_c(t)
 &=p(t)(I_k\otimes u_c)^*
       \Delta(t)^{-1}(I_k\otimes w_c)=0,\\
 \mathcal Q_c(t)
 &=p(t)(I_k\otimes w_c)^*
       \Delta(t)^{-1}(I_k\otimes u_c)=0.
 \end{aligned}
\]
These polynomial matrices vanish at infinitely many points, so
\[
 \mathcal P_c(t)\equiv0, \mathcal Q_c(t)\equiv0.
\]
In particular, $\mathcal P_c(0)=\mathcal Q_c(0)=0$.
Since $p(0)=\det\Delta\ne0$, we conclude that
\[
 \begin{aligned}
 (I_k\otimes u_c)^*\Delta^{-1}(I_k\otimes w_c)
 &=\frac{\mathcal P_c(0)}{p(0)}=0,\\
 (I_k\otimes w_c)^*\Delta^{-1}(I_k\otimes u_c)
 &=\frac{\mathcal Q_c(0)}{p(0)}=0.
 \end{aligned}
\]
This proves the claimed identities for every invertible $\Delta$.
\end{proof}

\subsection{Corners of Projections and Matrix Amplifications}

$\newline$ For $t\in\mathbb R$, define the Cayley transform of $2+t$ by
\begin{equation}\label{eq: 4.10}
 U(t)=\frac{2+t-i}{2+t+i},\qquad
 f(t)=\operatorname{Re}U(t)=\frac{(2+t)^2-1}{(2+t)^2+1}.
\end{equation}
We have $|U(t)|=1$. On $[0,1]$, $U(t)\notin\{1,-1\}$ and
$f'(t)=4(2+t)/((2+t)^2+1)^2>0$.

Let $x$ be a finite-spectrum positive contraction and write
$\sigma(x)=\{\lambda_1,\ldots,\lambda_m\}$. The Lagrange
interpolation formula gives
\[
 U(x)=\sum_{j=1}^m U(\lambda_j)
       \prod_{\substack{1\leq l\leq m\\ l\ne j}}
       \frac{x-\lambda_l I}{\lambda_j-\lambda_l}.
\]
To recover $x$ from the real part of $U(x)$, apply the Lagrange
interpolation formula at the distinct points $f(\lambda_j)$:
\[
 \psi_x(s)=\sum_{j=1}^m\lambda_j
   \prod_{\substack{1\leq l\leq m\\ l\ne j}}
     \frac{s-f(\lambda_l)}{f(\lambda_j)-f(\lambda_l)}.
\]
Since $\psi_x(f(\lambda))=\lambda$ for $\lambda\in\sigma(x)$,
finite-spectrum functional calculus gives
\begin{equation}\label{eq: 4.11}
\psi_x(f(x))=\psi_x(\Re U(x))=x.
\end{equation}

Let $p_0$ be a projection and let $h_0,k_0$ be finite-spectrum
positive contractions on a Hilbert space $\mathcal E$. Set
\[
 A_0=\mathcal A^*(I_{\mathcal E},p_0,h_0,k_0)\subseteq B(\mathcal E),
 \qquad x_1=p_0,\quad x_2=h_0,\quad x_3=k_0,
 \qquad U_c=U(x_c).
\]
We call $U_1,U_2,U_3$ the Cayley unitaries of the source tuple.
Lagrange interpolation and \eqref{eq: 4.11} give
$A_0=\mathcal A^*(I_{\mathcal E},U_1,U_2,U_3)$.

\begin{lemma}\label{lem: 4.4}
Define a linear map $v:\mathcal E\to \mathcal{K}\otimes \mathcal{E}$ by
\begin{equation}\label{eq:4.12}
 v\xi=\frac1{\sqrt6}\sum_{c=1}^3(u_c\otimes \xi +w_c\otimes U_c\xi),\qquad  (\xi\in\mathcal E),
\end{equation}
and set $P=vv^*$.
Then $v$ is an isometry, $P$ is a projection, and
${\mathcal A}^*(I_{\mathcal{K}\otimes \mathcal{E}},P,H\otimes I_{\mathcal E},K\otimes I_{\mathcal E})=M_R(A_0)$.
\end{lemma}

\begin{proof}
Since the six vectors $u_c,w_c$ are orthonormal and each $U_c$
is unitary, we have, for all $\xi,\eta\in\mathcal E$
\[
 \langle v\xi,v\eta\rangle
 =\frac16\sum_{c=1}^3
   \bigl(\langle\xi,\eta\rangle
         +\langle U_c\xi,U_c\eta\rangle\bigr)
 =\langle\xi,\eta\rangle.
\]
Hence $v^*v=I_{\mathcal E}$, and consequently
$P^*=P, P^2=vv^*vv^*=vv^*=P.$

Set
\(
 D=\mathcal A^*(I_{\mathcal K\otimes\mathcal E},
                    P,H\otimes I_{\mathcal E},K\otimes I_{\mathcal E}).
\)
Let $\{f_l\}_{l=1}^R$ be the orthonormal basis of $\mathcal K$
arising from the direct sum of the standard tensor product bases.
For each $l$, define the isometry
\[
 \widehat u_l:\mathcal E\longrightarrow
 \mathcal K\otimes\mathcal E,
 \qquad
 \widehat u_l\xi=f_l\otimes\xi.
\]
This choice of basis gives the $*$-isomorphism
\[
 B(\mathcal K\otimes\mathcal E)
 \longrightarrow M_R(B(\mathcal E)),
 \qquad
 T\longmapsto
 \bigl[\widehat u_l^*T\widehat u_k\bigr]_{l,k=1}^R.
\]
Under the identification, operators
on $\mathcal K\otimes\mathcal E$ are therefore represented by $R\times R$ operator
matrices.

Write
$u_c=\sum_{l=1}^R\alpha_{lc}f_l,~
 w_c=\sum_{l=1}^R\beta_{lc}f_l,~
 \alpha_{lc},\beta_{lc}\in\mathbb C.$
The definition of $v$ gives
\[
 \widehat u_l^*v
 =\frac1{\sqrt6}\sum_{c=1}^3
   \bigl(\alpha_{lc}I_{\mathcal E}+\beta_{lc}U_c\bigr)
 \in A_0.
\]
Thus $v\in M_{R,1}(A_0)$, that is to say,
$v:\mathcal E\to\mathcal K\otimes\mathcal E$
is represented by the column
\[
 v\longleftrightarrow
 \begin{pmatrix}
 \widehat u_1^*v\\
 \vdots\\
 \widehat u_R^*v
 \end{pmatrix}.
\]
In particular,
$\widehat u_l^*P\widehat u_k
 =(\widehat u_l^*v)(\widehat u_k^*v)^*
 \in A_0,$
so $P=vv^*\in M_R(A_0)$.
Since $H\otimes I_{\mathcal E}$ and
$K\otimes I_{\mathcal E}$ also belong to $M_R(A_0)$,
we have $D\subseteq M_R(A_0)$.

To prove the reverse inclusion, define
\[
 \theta:A_0\longrightarrow PM_R(A_0)P,
\qquad  \theta(a)=vav^*.
\]
The identity $v^*v=I_{\mathcal E}$ shows that $\theta$ is a
$*$-homomorphism with $\theta(1)=P$.
Lemma~\ref{lem:4.2} gives
$\mathcal A^*(1,H\otimes I_{\mathcal E},K\otimes I_{\mathcal E})
=B_0\otimes I_{\mathcal E}\subseteq D$.
For each $c$, let $b_c\in B_0$ act as
$E_{00}\otimes\Gamma\otimes I_N$ on the $c$th Hilbert space
summand of $\mathcal K$ and as zero on the other summands.
Then $b_cu_c=w_c,~ b_cw_c=u_c,~
 b_cu_d=b_cw_d=0~(d\ne c).$
The orthonormality of the vectors $u_c,w_c$ gives
$u_c^*b_cu_c=w_c^*b_cw_c=0$ and
$u_c^*b_cw_c=w_c^*b_cu_c=1$.
Therefore, with $\widetilde b_c:=b_c\otimes I_{\mathcal E}\in D$,
\[
 \begin{aligned}
 v^*\widetilde b_cv\xi
&= v^*(\frac1{\sqrt6}\sum_{d=1}^3
   \bigl(b_cu_d\otimes\xi+b_cw_d\otimes U_d\xi\bigr))\\
 &=v^*(\frac1{\sqrt6}
   \bigl(w_c\otimes\xi+u_c\otimes U_c\xi\bigr))\\
 &=\frac1{\sqrt6}\bigl(v^*(w_c\otimes\xi)
         +v^*(u_c\otimes U_c\xi)\bigr)\\
 &=\frac16\sum_{d=1}^3
   \bigl((u_d^*w_c)\xi+(w_d^*w_c)U_d^*\xi\bigr)\\
 &\quad+\frac16\sum_{d=1}^3
   \bigl((u_d^*u_c)U_c\xi+(w_d^*u_c)U_d^*U_c\xi\bigr)\\
 &=\frac16\bigl(U_c^*\xi+U_c\xi\bigr).
 \end{aligned}
\]
for every $\xi\in \mathcal E$.
By the arbitrariness of $\xi$, we conclude that
\[
 v^*\widetilde b_cv=\frac{U_c+U_c^*}{6}.
\]
Using $P=vv^*$ and $\theta(a)=vav^*$, we therefore have
\[
 \begin{aligned}
 3P\widetilde b_cP
 &=3vv^*\widetilde b_cvv^*\\
 &=3v\bigl(v^*\widetilde b_cv\bigr)v^*\\
 &=v\left(\frac{U_c+U_c^*}{2}\right)v^*\\
 &=v(\operatorname{Re}U_c)v^*\\
 &=\theta(\operatorname{Re}U_c).
 \end{aligned}
\]

Fix $c\in\{1,2,3\}$ and write
$\psi_{x_c}(t)=a_0+\sum_{j=1}^m a_jt^j$ with $a_j\in\mathbb C$.
By \eqref{eq: 4.11},
$x_c=a_0I_{\mathcal E}+\sum_{j=1}^m a_j(\operatorname{Re}U_c)^j$.
Since $\theta$ is a $*$-homomorphism and
$\theta(I_{\mathcal E})=P$, we obtain
\[
 \begin{aligned}
 \theta(x_c)
 &=a_0P+\sum_{j=1}^{m}
   a_j\theta\bigl((\operatorname{Re}U_c)^j\bigr)\\
 &=a_0P+\sum_{j=1}^{m}
   a_j\bigl(\theta(\operatorname{Re}U_c)\bigr)^j\\
 &=a_0P+\sum_{j=1}^{m}
   a_j(3P\widetilde b_cP)^j
 \in D.
 \end{aligned}
\]

Since $A_0=\mathcal A^*(I_{\mathcal E},x_1,x_2,x_3)$ and the
$x_c$ are selfadjoint, each $a\in A_0$ has the form
\[
 a=\lambda_0I_{\mathcal E}
   +\sum_{\nu=1}^s\lambda_\nu x_{i_{\nu,1}}\cdots x_{i_{\nu,r_\nu}},
 \qquad \lambda_\nu\in\mathbb C,\quad i_{\nu,j}\in\{1,2,3\}.
\]
It follows that
\[
 \theta(a)=\lambda_0P+
 \sum_{\nu=1}^s\lambda_\nu\theta(x_{i_{\nu,1}})\cdots
                              \theta(x_{i_{\nu,r_\nu}})\in D.
\]
Therefore $\theta(A_0)\subseteq D$.

Set
\[
 \mathcal W=\{dv:d\in D\}
 \subseteq B(\mathcal E,\mathcal K\otimes\mathcal E).
\]
For every $d\in D$ and $a\in A_0$,
$(dv)a=dvav^*v=d\theta(a)v\in\mathcal W$. Thus $\mathcal W$
is a right $A_0$-module.

We claim that $\widehat u_l\in\mathcal W$ for every $l$.
Write
\[
 \Gamma=\operatorname{diag}(\varepsilon_1,\ldots,\varepsilon_N),
 \qquad  \varepsilon_b\in\{1,-1\}.
\]
Fix $l$. There are indices $c,t,\alpha,b$ such that $f_l$
has component $e_t\otimes e_\alpha\otimes e_b$ in the $c$th
Hilbert space summand of $\mathcal K$ and zero components
elsewhere. Let $F_l\in B_0$ act as
$E_{t0}\otimes E_{\alpha b}\otimes I_N$
on that summand and as zero on the others. By
the definitions of $\Omega,u_c,w_c$, we have
\[
 F_lu_c=\frac1{\sqrt N}f_l,\qquad
 F_lw_c=\frac{\varepsilon_b}{\sqrt N}f_l,\qquad
 F_lu_d=F_lw_d=0\quad(d\ne c).
\]
Consequently,
\[
 (F_l\otimes I_{\mathcal E})v
 =\frac1{\sqrt{6N}}\,
   \widehat u_l(I_{\mathcal E}+\varepsilon_bU_c).
\]

Both $I_{\mathcal E}-U_c$ and $I_{\mathcal E}+U_c$
have inverses in $A_0$. Set $T_c=2I_{\mathcal E}+x_c$. Both
$T_c$ and $T_c+iI_{\mathcal E}$ are invertible. Since $x_c$ has
finite spectrum, there is a polynomial $r_c$ satisfying
$r_c(\lambda)=(2+\lambda)^{-1}$ for every
$\lambda\in\sigma(x_c)$. Hence
$T_c^{-1}=r_c(x_c)\in A_0$.
Since $U_c=(T_c-iI_{\mathcal E})(T_c+iI_{\mathcal E})^{-1}$, it follows that
\[
\begin{aligned}
I_{\mathcal E}-U_c
&=(T_c+iI_{\mathcal E})(T_c+iI_{\mathcal E})^{-1}
  -(T_c-iI_{\mathcal E})(T_c+iI_{\mathcal E})^{-1}\\
&=\bigl[(T_c+iI_{\mathcal E})-(T_c-iI_{\mathcal E})\bigr](T_c+iI_{\mathcal E})^{-1}\\
&=2i(T_c+iI_{\mathcal E})^{-1}.
\end{aligned}
\]
Hence
$(I_{\mathcal E}-U_c)^{-1}
 =\frac{T_c+iI_{\mathcal E}}{2i}\in A_0$, and similarly $(I_{\mathcal E}+U_c)^{-1}
 =\frac12\bigl(I_{\mathcal E}+iT_c^{-1}\bigr)\in A_0$.
The right module property of $\mathcal W$ also implies
\[
 \widehat u_l
 =\sqrt{6N}\,(F_l\otimes I_{\mathcal E})v
   (I_{\mathcal E}+\varepsilon_bU_c)^{-1}
 \in\mathcal W.
\]
This proves the claim.

For arbitrary $a\in A_0$ and $1\leq l,k\leq R$,
we now have $\widehat u_l a,\widehat u_k\in\mathcal W$.
Choose $d,e\in D$ such that
$\widehat u_l a=dv$ and $\widehat u_k=ev.$
Then $\widehat u_l a\widehat u_k^*=dvv^*e^*
 =dPe^*\in D.$
Under the matrix identification above, this operator has $a$
in its $(l,k)$ entry and zero in every other entry, since
$\widehat u_i^*
 (\widehat u_l a\widehat u_k^*)
 \widehat u_j
 =\delta_{il}\delta_{kj}a.$
These elementary operator matrices span $M_R(A_0)$.
Thus $M_R(A_0)\subseteq D$, and hence
\[
 {\mathcal A}^*(I_{\mathcal{K}\otimes \mathcal{E}},P,H\otimes I_{\mathcal E},K\otimes I_{\mathcal E})=M_R(A_0).
\]
\end{proof}

\begin{corollary}\label{cor 4.5}
Let $k\geq1$ and suppose
$\Delta=Z_0\otimes I_R-X_0\otimes H-Y_0\otimes K$ is invertible,
where $Z_0,X_0,Y_0\in M_k(\mathbb C)$. Put $V=I_k\otimes v$,
$L_c=I_k\otimes u_c$, and $M_c=I_k\otimes w_c$. Then
\[
 \begin{aligned}
 V^*(\Delta^{-1}\otimes I_{\mathcal E})V
 =\frac16\sum_{c=1}^3\Bigl(
 &(L_c^*\Delta^{-1}L_c+M_c^*\Delta^{-1}M_c)\otimes I_{\mathcal E}\\
 &+(L_c^*\Delta^{-1}M_c)\otimes U_c
  +(M_c^*\Delta^{-1}L_c)\otimes U_c^*
 \Bigr).
 \end{aligned}
\]
In particular, if $k\leq n$, then
\[
 V^*(\Delta^{-1}\otimes I_{\mathcal E})V
 =\frac13\sum_{c=1}^3(L_c^*\Delta^{-1}L_c)\otimes I_{\mathcal E}.
\]
\end{corollary}

\begin{proof}
By \eqref{eq:4.12},
\[
 V=I_k\otimes v
 =\frac1{\sqrt6}\sum_{c=1}^3
 \bigl(L_c\otimes I_{\mathcal E}+M_c\otimes U_c\bigr),
\]
and likewise
\[
 V^*
 =\frac1{\sqrt6}\sum_{c=1}^3
 \bigl(L_c^*\otimes I_{\mathcal E}+M_c^*\otimes U_c^*\bigr).
\]
Consequently,
\[
 \begin{aligned}
 V^*(\Delta^{-1}\otimes I_{\mathcal E})V
 &=\frac16\sum_{c,d=1}^3
 \Bigl[
 (L_c^*\Delta^{-1}L_d)\otimes I_{\mathcal E}
  +(L_c^*\Delta^{-1}M_d)\otimes U_d\\
  &\qquad +(M_c^*\Delta^{-1}L_d)\otimes U_c^*
       +(M_c^*\Delta^{-1}M_d)\otimes U_c^*U_d
 \Bigr].
 \end{aligned}
\]

Since $\Delta^{-1}$ preserves the summands of
$\mathbb C^k\otimes\mathcal K$ and $u_c,w_c$ are supported in
$\mathcal K_c$, we have, for $c\ne d$,
\begin{equation}\label{eq: 4.13}
 \begin{aligned}
 L_c^*\Delta^{-1}L_d&=0,\qquad L_c^*\Delta^{-1}M_d=0,\\
 M_c^*\Delta^{-1}L_d&=0,\qquad M_c^*\Delta^{-1}M_d=0.
 \end{aligned}
\end{equation}
Hence \eqref{eq: 4.13} and $U_c^*U_c=I_{\mathcal E}$ give
\[
\begin{aligned}
 V^*(\Delta^{-1}\otimes I_{\mathcal E})V
 &=\frac16\sum_{c=1}^3
  \Bigl(
  (L_c^*\Delta^{-1}L_c)\otimes I_{\mathcal E}
  +(L_c^*\Delta^{-1}M_c)\otimes U_c\\
 &\qquad
  +(M_c^*\Delta^{-1}L_c)\otimes U_c^*
  +(M_c^*\Delta^{-1}M_c)\otimes U_c^*U_c
  \Bigr)\\
 &=\frac16\sum_{c=1}^3
  \Bigl(
  (L_c^*\Delta^{-1}L_c+M_c^*\Delta^{-1}M_c)\otimes I_{\mathcal E}\\
 &\qquad
  +(L_c^*\Delta^{-1}M_c)\otimes U_c
  +(M_c^*\Delta^{-1}L_c)\otimes U_c^*
  \Bigr).
 \end{aligned}
\]

If $k\leq n$, then Lemma~\ref{lem: 4.3} and \eqref{eq:4.9}
give
\[
 V^*(\Delta^{-1}\otimes I_{\mathcal E})V
 =\frac13\sum_{c=1}^3 L_c^*\Delta^{-1}L_c\otimes I_{\mathcal E}.
\]
\end{proof}

\subsection{Iterated Amplifications and Finite-Dimensional Representations}

\(\newline\)For each integer $j\geq1$, set
\begin{equation}
 \kappa_j=j14^{j-1},~
 g_j=\min\{g\in2\mathbb N:g>\kappa_j^2\},~
 N_j=2^{g_j/2},~
 R_j=3(2g_j+1)N_j^2.
\end{equation}
We apply the preceding construction with $n=\kappa_j$ and denote its
Hilbert space by
\[
 \mathcal K_j
 =\bigoplus_{c=1}^3
 \bigl({\mathbb{C}}^{V_{g_j}}\otimes{\mathbb{C}}^{N_j}\otimes{\mathbb{C}}^{N_j}\bigr),
  \qquad\dim\mathcal K_j=R_j.
\]
Write $H_j,K_j\in B(\mathcal K_j)$ and
$u_{j,c},w_{j,c}\in\mathcal K_j$ for the corresponding operators
and vectors, with $1\leq c\leq3$.
Then $H_j$ and $K_j$ are fixed finite-dimensional positive
contractions, and the six vectors $u_{j,c},w_{j,c}$ are
orthonormal.

Call a tuple $x=(x_1,x_2,x_3)$ on a Hilbert space $\mathcal E$
admissible if $x_1$ is a projection and $x_2,x_3$ are
finite-spectrum positive contractions. For an admissible tuple, define
\[
 v_j(x):\mathcal E\longrightarrow\mathcal K_j\otimes\mathcal E
\]
by
\[
 v_j(x)\xi
 =\frac1{\sqrt6}\sum_{c=1}^3
 \bigl(u_{j,c}\otimes\xi+w_{j,c}\otimes U(x_c)\xi\bigr),
 \qquad \xi\in\mathcal E,
\]
and set
\begin{equation}\label{eq:encoder}
 E_j(x)=
 \bigl(v_j(x)v_j(x)^*,\,
       H_j\otimes I_{\mathcal E},\,
       K_j\otimes I_{\mathcal E}\bigr).
\end{equation}
By Lemma~\ref{lem: 4.4}, the first component is a projection. The other two
components are amplifications of finite-dimensional positive
contractions and hence remain positive contractions with
finite spectra. Thus $E_j(x)$ is an admissible tuple on
$\mathcal K_j\otimes\mathcal E$.

Fix $m\geq 1$. We construct an admissible tuple
$T_j^{(m)}$ on a finite-dimensional Hilbert space
$\mathcal E_j^{(m)}$ by descending induction on $j$.
Write $D_j^{(m)}=\dim\mathcal E_j^{(m)}$.

At the base stage $j=m+1$, set
\[
 \mathcal E_{m+1}^{(m)}=\mathbb C,~
 T_{m+1}^{(m)}=(0,0,0),\qquad
 D_{m+1}^{(m)}=1.
\]
Each zero denotes the zero operator on $\mathbb C$.
At the first step $j=m$, \eqref{eq: 4.10} gives
\[
 U(0)=\frac{2-i}{2+i}=\frac{3-4i}{5}.
\]
Under the natural identification
$\mathcal K_m\otimes\mathbb C\cong\mathcal K_m$, the isometry
$v_m(T_{m+1}^{(m)})$ sends $\xi\in\mathbb C$ to $\xi z_m$, where
\[
 z_m=\frac1{\sqrt6}\sum_{c=1}^3
 \left(u_{m,c}+\frac{3-4i}{5}\,w_{m,c}\right).
\]
The six vectors $u_{m,c},w_{m,c}$ are orthonormal and
$|U(0)|=1$, so $\|z_m\|=1$. Hence, under the same identification,
\[
 T_m^{(m)}
 =E_m\bigl(T_{m+1}^{(m)}\bigr)
 =\bigl(z_mz_m^*,H_m,K_m\bigr).
\]

Here $z_mz_m^*$ is the rank-one orthogonal projection onto
$\mathbb C z_m$. Under the canonical identification
\(\mathcal K_m\otimes\mathbb C\cong\mathcal K_m\), the operators
\(H_m\otimes I_{\mathbb C}\) and \(K_m\otimes I_{\mathbb C}\)
correspond to \(H_m\) and \(K_m\), respectively.

In general, given the admissible tuple $T_{j+1}^{(m)}$ on
$\mathcal E_{j+1}^{(m)}$, define
\[
 \mathcal E_j^{(m)}
 =\mathcal K_j\otimes\mathcal E_{j+1}^{(m)},~
 T_j^{(m)}=E_j\bigl(T_{j+1}^{(m)}\bigr).
\]

At each step, the dimension satisfies
$D_j^{(m)}=R_jD_{j+1}^{(m)}.$
Since $D_{m+1}^{(m)}=1$, we have
$D_j^{(m)}=\prod_{i=j}^m R_i\quad (1\leq j\leq m)$.
Write
\[
 T_1^{(m)}=(P^{(m)},H^{(m)},K^{(m)}),
~D_m=D_1^{(m)}.
\]
Since $\mathcal E_1^{(m)}=\mathcal K_1\otimes\mathcal E_2^{(m)}$, $
D_m=D_1^{(m)}=R_1D_2^{(m)}$, and so
$\dim\mathcal E_2^{(m)}=D_2^{(m)}=\frac{D_m}{R_1}.$ Therefore, in the tensor product basis,
\begin{equation}\label{eq: 4.16}
 H^{(m)}=H_1\otimes I_{D_m/R_1},~
 K^{(m)}=K_1\otimes I_{D_m/R_1}.
\end{equation}

\begin{lemma}\label{lem: 4.6}
For every $m\geq1$, the tensor product basis identifies
$B(\mathcal E_1^{(m)})$ with $M_{D_m}(\mathbb C)$. Under this
identification,
\[
 \mathcal A^*\bigl(
 I_{\mathcal E_1^{(m)}},P^{(m)},H^{(m)},K^{(m)}
 \bigr)
 =M_{D_m}(\mathbb C).
\]
\end{lemma}

\begin{proof}
Fix $m$. For $1\leq j\leq m+1$, let $A_j^{(m)}$ be the unital
$*$-algebra generated by the three components of $T_j^{(m)}$
on $\mathcal E_j^{(m)}$.

At $j=m+1$, the zero tuple acts on $\mathbb C$, so
$A_{m+1}^{(m)}=\mathbb C$. Suppose
$A_{j+1}^{(m)}=M_{D_{j+1}^{(m)}}(\mathbb C)$.
By Lemma~\ref{lem: 4.4} and $D_j^{(m)}=R_jD_{j+1}^{(m)}$, we have
\[
 A_j^{(m)}
 =M_{R_j}\bigl(A_{j+1}^{(m)}\bigr)
 =M_{R_j}\bigl(M_{D_{j+1}^{(m)}}(\mathbb C)\bigr)
 =M_{D_j^{(m)}}(\mathbb C).
\]
Taking $j=1$ gives the claim.
\end{proof}

\begin{remark}\label{rem: 4.7}
For each \(m\geq 1\), we start from the same one-dimensional
base data,
\(\mathcal E_{m+1}^{(m)}=\mathbb C,
~ T_{m+1}^{(m)}=(0,0,0),
 ~D_{m+1}^{(m)}=1.\)
Applying \(E_m,E_{m-1},\ldots,E_1\) in that order gives
\[
 T_1^{(m)}=(E_1\circ\cdots\circ E_m)(0,0,0).
\]
Each step uses the triangular bouquets at its stage and adds
one layer to the construction. Thus \(T_1^{(m)}\) is obtained
from the same base tuple by a construction of depth \(m\).

For $m=2$, the base stage gives $T_2^{(2)}=E_2(0,0,0)$
on $\mathcal K_2$. Thus the input to $E_1$ has dimension $R_2$,
and the output tuple
\[
 T_1^{(2)}
 =E_1\bigl(T_2^{(2)}\bigr)
 =\left(
 P^{(2)},\
 H_1\otimes I_{R_2},\
 K_1\otimes I_{R_2}
 \right)
\]
acts on $\mathcal K_1\otimes\mathcal K_2$.

For $m=3$, we have $T_3^{(3)}=E_3(0,0,0)$ on $\mathcal K_3$.
Applying $E_2$ then yields
\[
 T_2^{(3)}
 =E_2\bigl(T_3^{(3)}\bigr)
 =\left(
 P_2^{(3)},\
 H_2\otimes I_{R_3},\
 K_2\otimes I_{R_3}
 \right)
\]
on $\mathcal K_2\otimes\mathcal K_3$.
The input to $E_1$ now acts on $\mathcal K_2\otimes\mathcal K_3$,
of dimension $R_2R_3$. Therefore
\[
 T_1^{(3)}
 =E_1\bigl(T_2^{(3)}\bigr)
 =\left(
 P^{(3)},\
 H_1\otimes I_{R_2R_3},\
 K_1\otimes I_{R_2R_3}
 \right)
\]
acts on
$\mathcal K_1\otimes\mathcal K_2\otimes\mathcal K_3$.
\end{remark}

On $\mathcal H=\bigoplus_{m\geq1}\mathbb C^{D_m}$, define
\begin{equation}\label{eq: 4.17}
 p=\bigoplus_mP^{(m)},\qquad h=\bigoplus_mH^{(m)},\qquad
 k=\bigoplus_mK^{(m)},\qquad A=\mathcal A^*(I_{\mathcal H},p,h,k).
\end{equation}
These operators are positive contractions, and $p$ is a projection.
By \eqref{eq: 4.16}, $h$ and $k$ are amplifications of the fixed
finite matrices $H_1$ and $K_1$, so all three generators have finite
spectra. Lemma~\ref{lem: 4.6} shows that coordinate evaluation maps
$A$ onto $M_{D_m}(\mathbb C)$. Since
$R_j\geq 3\cdot5\cdot2^2=60$ and $D_m\geq60^m\to\infty$, the
algebra $A$ is infinite dimensional. To prove its regularity, we
study the spectra of selfadjoint matrix polynomials in its generators.

Let $\mathcal P=\mathbb C\langle x_1,x_2,x_3\rangle$ be the
free unital algebra of noncommutative polynomials in the
formal symbols $x_1,x_2,x_3$, equipped with
the involution defined by $x_r^*=x_r$ for $r=1,2,3$.
The algebra $M_\ell(\mathcal P)$ consists of $\ell\times\ell$
matrices with entries in $\mathcal P$, with involution given by $[a_{ij}]^*=[a_{ji}^*]$.

\subsection{Linearization and Schur complements}

\(\newline\)We use selfadjoint linearization by Schur complements
\cite[Propositions~3 and~4]{Anderson}.

\begin{lemma}\label{lem: 4.8}
Let $a\in M_\ell(\mathcal P)$ satisfy $a^*=a$, and put
$a_0=a(0,0,0)$. There exist an integer $n_1\geq\ell$, a constant matrix
$V\in M_{\ell,n_1-\ell}(\mathbb C)$, and a selfadjoint
matrix polynomial $B(x)\in M_{n_1-\ell}(\mathcal P)$
whose entries have degree at most one, such that
$B(x)^{-1}\in M_{n_1-\ell}(\mathcal P)$,
and the matrix polynomial
\[
 L_1(z;x)=
 \begin{pmatrix}
 zI_\ell-a_0&V\\
 V^*&B(x)
 \end{pmatrix},
 ~ z\in\mathbb C,
\]
has Schur complement
$zI_\ell-a_0-VB(x)^{-1}V^*
 =zI_\ell-a(x)$.
Moreover,
$L_1(z;x)$ can be written as
 $L_1(z;x)
 =C_0+zC_z+C_1x_1+C_2x_2+C_3x_3
 \in M_{n_1}(\mathcal P)$ where
$C_0,C_z,C_1,C_2,C_3\in M_{n_1}(\mathbb C)$ are selfadjoint matrices. For constant $a$, one may take
$n_1=\ell$ and omit the lower and off-diagonal blocks.
\end{lemma}
\begin{proof}
If $a$ is constant, take $n_1=\ell$ and
$L_1(z;x)=zI_\ell-a$.
Assume now that $a$ is nonconstant.
Write
\[
 a(x)=a_0+\sum_{w\in\mathcal W}A_wx_w,~ a_0=a(0,0,0),\quad A_w\in M_\ell(\mathbb C),
\]
where $\mathcal W$ is a finite set of nonempty words.
For $w=(i_1,\ldots,i_d)$, we write
$x_w=x_{i_1}\cdots x_{i_d}$.
Since $a^*=a$, we have $a_0^*=a_0$.

For each word $w=(i_1,\ldots,i_d)$, define a block matrix
$N_w(x)$ with block indices $0,\ldots,d$ and entries
\[
 (N_w(x))_{rs}
 =
 \begin{cases}
 x_{i_{r+1}}I_\ell,&s=r+1,\\
 0,&\text{otherwise}.
 \end{cases}
\]
This matrix is strictly upper triangular, so
$N_w(x)^{d+1}=0$. Consequently,
$T_w(x):=I_{(d+1)\ell}-N_w(x)$
has the polynomial inverse
$T_w(x)^{-1}
 =\sum_{r=0}^{d}N_w(x)^r.
$
Thus its $(0,d)$ block is
\[
 \bigl(T_w(x)^{-1}\bigr)_{0d}
 =(x_{i_1}I_\ell)\cdots(x_{i_d}I_\ell)
 =x_wI_\ell.
\]

Set
\[
 b_w=\begin{pmatrix}A_w&0&\cdots&0\end{pmatrix},~~
c_w=
\begin{pmatrix}
0\\
\vdots\\
0\\
I_\ell
\end{pmatrix}.
\]
Then
$b_wT_w(x)^{-1}c_w=A_wx_w.$
Choose an ordering $\mathcal W=\{w_1,\ldots,w_t\}$ and set
\[
 \begin{gathered}
 q=\ell\sum_{w\in\mathcal W}(|w|+1),~
 T(x)=\operatorname{diag}
       \bigl(T_{w_1}(x),\ldots,T_{w_t}(x)\bigr),\\
 b=\begin{pmatrix}b_{w_1}&\cdots&b_{w_t}\end{pmatrix},
 ~
c=
\begin{pmatrix}
c_{w_1}\\
\vdots\\
c_{w_t}
\end{pmatrix}
 \end{gathered}
\]
Thus $T(x),T(x)^{-1}\in M_q(\mathcal P)$ and
$bT(x)^{-1}c=a(x)-a_0.$

Now put
\[
 B(x)=
 \begin{pmatrix}
 0&T(x)\\
 T(x)^*&0
 \end{pmatrix},
 \qquad
 V=\frac1{\sqrt2}\begin{pmatrix}c^*&b\end{pmatrix}.
\]
It is clear that
\[
 B(x)^{-1}
 =
 \begin{pmatrix}
 0&(T(x)^*)^{-1}\\
 T(x)^{-1}&0
 \end{pmatrix}
 \in M_{2q}(\mathcal P).
\]
Moreover,
\[
 \begin{aligned}
 VB(x)^{-1}V^*
 &=\frac12\left(
 c^*(T(x)^*)^{-1}b^*+bT(x)^{-1}c
 \right)\\
 &=\frac12\left(
 (a(x)-a_0)^*+a(x)-a_0
 \right)\\
 &=a(x)-a_0.
 \end{aligned}
\]
Set $n_1=\ell+2q$ and define
\begin{equation}\label{eq: 4.18}
 L_1(z;x)
 =
 \begin{pmatrix}
 zI_\ell-a_0&V\\
 V^*&B(x)
 \end{pmatrix}.
\end{equation}
The preceding identity yields
$zI_\ell-a_0-VB(x)^{-1}V^*
 =zI_\ell-a(x).$

We now identify the scalar coefficient matrices of $L_1$.
Set
\[
(S_{w,r})_{uv}
=
\begin{cases}
I_\ell,&v=u+1,\ i_{u+1}=r,\\
0,&\text{otherwise},
\end{cases}
\]
and write
$S_r=\operatorname{diag}
\bigl(S_{w_1,r},\ldots,S_{w_t,r}\bigr)\in M_q(\mathbb C).$
Then
\[
\begin{aligned}
T(x)
&=\operatorname{diag}(T_{w_1}(x),\ldots,T_{w_t}(x))\\
&=I_q-\sum_{r=1}^3S_rx_r.
\end{aligned}
\]
Substituting the expressions for $B(x)$ and $V$ into \eqref{eq: 4.18} gives
\[
L_1(z;x)=
\begin{pmatrix}
zI_\ell-a_0&c^*/\sqrt2&b/\sqrt2\\
c/\sqrt2&0&I_q-\sum_{r=1}^3S_rx_r\\
b^*/\sqrt2&I_q-\sum_{r=1}^3S_r^*x_r&0
\end{pmatrix}.
\]
Therefore
$L_1(z;x)=C_0+zC_z+\sum_{r=1}^3C_rx_r,$
where
\[
 C_0=
 \begin{pmatrix}
 -a_0&c^*/\sqrt2&b/\sqrt2\\
 c/\sqrt2&0&I_q\\
 b^*/\sqrt2&I_q&0
 \end{pmatrix},
 \qquad
 C_z=\operatorname{diag}(I_\ell,0_q,0_q),
\]
and
\[
 C_r=
 \begin{pmatrix}
 0&0&0\\
 0&0&-S_r\\
 0&-S_r^*&0
 \end{pmatrix}
 \qquad(1\leq r\leq3).
\]
Since all these coefficient matrices are selfadjoint, it follows that
$L_1(z;x)^*=L_1(\overline z;x).$
\end{proof}

\begin{remark}\label{rem: 4.9}
For nonconstant $a$, Lemma~\ref{lem: 4.8} gives a polynomial inverse
of $B(x)$. Block elimination gives
\[
 \begin{pmatrix}
 I_\ell&-VB(x)^{-1}\\
 0&I_{2q}
 \end{pmatrix}
 L_1(z;x)
 \begin{pmatrix}
 I_\ell&0\\
 -B(x)^{-1}V^*&I_{2q}
 \end{pmatrix}
 =
 \begin{pmatrix}
 zI_\ell-a(x)&0\\
 0&B(x)
 \end{pmatrix},
\]
The two triangular factors are invertible. Thus $L_1(z;x)$ is
invertible if and only if $zI_\ell-a(x)$ is invertible.
For bounded selfadjoint operators $x_1,x_2,x_3$ on a Hilbert space
$\mathcal E$, evaluation gives
\[
 \begin{aligned}
 L_1(z;x)
 &=
 \begin{pmatrix}
 (zI_\ell-a_0)\otimes I_{\mathcal E}&V\otimes I_{\mathcal E}\\
 V^*\otimes I_{\mathcal E}&B(x)
 \end{pmatrix}\\
 &=C_0\otimes I_{\mathcal E}+zC_z\otimes I_{\mathcal E}
   +\sum_{r=1}^3 C_r\otimes x_r
 \in M_{n_1}(B(\mathcal E)),
 \qquad z\in\mathbb C,
 \end{aligned}
\]
where
\[
 B(x)=
 \begin{pmatrix}
 0&I_q\otimes I_{\mathcal E}-\sum_{r=1}^3 S_r\otimes x_r\\
 I_q\otimes I_{\mathcal E}-\sum_{r=1}^3 S_r^*\otimes x_r&0
 \end{pmatrix}.
\]
Since evaluation preserves products, evaluating the polynomial
inverse \(B^{-1}\) at \(x\) gives the inverse of the lower block
\(B(x)\). The elimination identity therefore implies that
\(L_1(z;x)\) is invertible if and only if
\(zI_{\mathbb C^\ell\otimes\mathcal E}-a(x)\) is invertible.
For constant $a$, the same equivalence follows directly from the
formula for $L_1$ with the lower block omitted.
\end{remark}

\begin{lemma}\label{lem: 4.10}
Let $\Delta\in M_b(B(\mathcal E))$ be invertible,
$V\in M_{b,q}(B(\mathcal E))$, and $B=B^*\in M_q(\mathbb{C})$.
Put $C=V^*\Delta^{-1}V$. Then $\Delta-VBV^*$ is invertible if and only if
$F=\begin{pmatrix}-C&I_q\\I_q&-B\end{pmatrix}$
is invertible. Here scalar matrices are amplified by $I_{\mathcal E}$.
\end{lemma}

\begin{proof}
Let
\[
 J=\begin{pmatrix}0&I_q\\I_q&-B\end{pmatrix},\qquad
 J^{-1}=\begin{pmatrix}B&I_q\\I_q&0\end{pmatrix},\qquad
 W=[\,V,0\,],
\]
and consider $\mathcal G=\begin{pmatrix}\Delta&W\\W^*&J\end{pmatrix}$.
Taking the Schur complement of $J$ gives
$\Delta-WJ^{-1}W^*=\Delta-VBV^*$.
Taking instead the Schur complement of $\Delta$ gives
\[
 J-W^*\Delta^{-1}W
 =\begin{pmatrix}-C&I_q\\I_q&-B\end{pmatrix}=F.
\]
Since $J$ and $\Delta$ are invertible, the two Schur complement
identities give
\[
 \begin{aligned}
 \Delta-VBV^*\text{ is invertible}
 &\ \Longleftrightarrow\ \mathcal G\text{ is invertible}\\
 &\ \Longleftrightarrow\
 J-W^*\Delta^{-1}W=F\text{ is invertible}.
 \end{aligned}
\]
\end{proof}

\begin{lemma}\label{lem: 4.11}
Let $x_1,x_2,x_3$ be bounded selfadjoint operators on a Hilbert
space $\mathcal E$. Use the Cayley transform \eqref{eq: 4.10}
on bounded selfadjoint operators. For
$C_0,C_1,C_2,C_3\in M_q(\mathbb C)$, set
\[
 F=C_0\otimes I_{\mathcal{E}}+\sum_{r=1}^3
       \bigl(C_r\otimes U(x_r)+C_r^*\otimes U(x_r)^*\bigr),
 \qquad C_0=C_0^*.
\]
Then \(F\) is the Schur complement of the lower-right block of
a selfadjoint \(7q\times 7q\) matrix polynomial of degree at most
one in \(x_1,x_2,x_3\). The lower-right block is invertible at
every bounded selfadjoint substitution.
\end{lemma}

\begin{proof}
Put
$D_0=C_0+\sum_{r=1}^3(C_r+C_r^*),~
 E_r=\begin{pmatrix}I_q&-2iC_r^*\end{pmatrix},
~\widetilde E_r=E_r\otimes I_{\mathcal E}.$
Set $T_r=2I_{\mathcal E}+x_r$ and define
\[
 M_r(x)=
 \begin{pmatrix}
 0&I_q\otimes(T_r-iI_{\mathcal E})\\
 I_q\otimes(T_r+iI_{\mathcal E})&0
 \end{pmatrix}.
\]
Since $T_r$ is selfadjoint, $T_r\pm iI_{\mathcal E}$
are invertible, and hence
\[
 M_r(x)^{-1}=
 \begin{pmatrix}
 0&I_q\otimes(T_r+iI_{\mathcal E})^{-1}\\
 I_q\otimes(T_r-iI_{\mathcal E})^{-1}&0
 \end{pmatrix}.
\]

Define
\[
 T(x)=
 \begin{pmatrix}
 D_0\otimes I_{\mathcal E}
   &\widetilde E_1&\widetilde E_2&\widetilde E_3\\
 \widetilde E_1^*&M_1(x)&0&0\\
 \widetilde E_2^*&0&M_2(x)&0\\
 \widetilde E_3^*&0&0&M_3(x)
 \end{pmatrix}.
\]
The Cayley transform gives
\[
 (T_r+iI_{\mathcal E})^{-1}=\frac{I_{\mathcal E}-U(x_r)}{2i},
 \qquad
 (T_r-iI_{\mathcal E})^{-1}=\frac{U(x_r)^*-I_{\mathcal E}}{2i}.
\]
Consequently,
\[
 \begin{aligned}
 \widetilde E_rM_r(x)^{-1}\widetilde E_r^*
 &=
 2iC_r\otimes(T_r+iI_{\mathcal E})^{-1}
 -2iC_r^*\otimes(T_r-iI_{\mathcal E})^{-1}\\
 &=
 (C_r+C_r^*)\otimes I_{\mathcal E}
 -C_r\otimes U(x_r)-C_r^*\otimes U(x_r)^*.
 \end{aligned}
\]
Therefore
$D_0\otimes I_{\mathcal E}
 -\sum_{r=1}^3
 \widetilde E_rM_r(x)^{-1}\widetilde E_r^*
 =F(x),$
which is the Schur complement of the invertible lower block
$\operatorname{diag}(M_1(x),M_2(x),M_3(x))$ in $T(x)$.
The displayed blocks also give selfadjointness and auxiliary size
$q+3(2q)=7q$.
\end{proof}

\subsection{Uniform Spectral Finiteness and Three-Projection Generation}

\(\newline\)Let $\mathbb C(z)$ denote the field of rational functions
in the scalar variable $z$:
\[
 \mathbb C(z)
 =\left\{\frac{p(z)}{q(z)}:p,q\in\mathbb C[z],\ q\not\equiv0\right\}.
\]
For $n\geq1$, $M_n(\mathbb C(z))$ consists of the $n\times n$
matrices with entries in $\mathbb C(z)$. For $R(z)=[r_{ij}(z)]\in M_n(\mathbb C(z))$, define
\[
 R^\#(z)=R(\overline z)^*
 =\bigl[\overline{r_{ji}(\overline z)}\bigr]_{i,j}.
\]
Thus $R^\#$ is again rational, and $R^\#(t)=R(t)^*$ at every
real $t$ where $R$ is defined.

Fix $a=a^*\in M_\ell(\mathcal P)$ and let $x=(x_1,x_2,x_3)$ act
on $\mathcal E$. We construct matrix polynomials of the form
\begin{equation}\label{eq: 4.19}
 \begin{aligned}
 L_j(z;x)
 ={}&Z_j(z)\otimes I_{\mathcal E}+A_j(z)\otimes x_1
 +B_j(z)\otimes x_2+C_j(z)\otimes x_3,
 \end{aligned}
\end{equation}
with $Z_j,A_j,B_j,C_j\in M_{n_j}(\mathbb C(z))$ satisfying
\[
 Z_j^\#=Z_j,\qquad A_j^\#=A_j,\qquad
 B_j^\#=B_j,\qquad C_j^\#=C_j.
\]
At any $z$ away from the coefficient poles, $L_j(z;x)$ belongs
to $M_{n_j}(B(\mathcal E))$ and has degree at most one in the
source variables. Its coefficients depend only on \(a\) and on the matrices \(H_i,K_i\) and
compression vectors \(u_{i,c},w_{i,c}\) with
\(1\leq i<j\) and \(1\leq c\leq 3\). We call these matrices and vectors the first \(j-1\) backgrounds.

Write $C_r^{\mathrm{lin}}=C_r$ for the constant coefficients
of Lemma~\ref{lem: 4.8}, to distinguish them from the functions
$C_j(z)$. By Remark~\ref{rem: 4.9},
\[
 \begin{aligned}
 L_1(z;x)
 &=Z_1(z)\otimes I_{\mathcal E}+A_1(z)\otimes x_1
   +B_1(z)\otimes x_2+C_1(z)\otimes x_3,\\
 Z_1(z)&=C_0+zC_z,\qquad A_1(z)=C_1^{\mathrm{lin}},\\
 B_1(z)&=C_2^{\mathrm{lin}},\qquad C_1(z)=C_3^{\mathrm{lin}}.
 \end{aligned}
\]
Its coefficients depend only on $a$. Once $L_j$ has been
constructed, define
\begin{equation}\label{eq: 4.20}
 \begin{aligned}
 \Delta_j(z)&=L_j(z;0,H_j,K_j)\\
 &=Z_j(z)\otimes I_{\mathcal K_j}
   +B_j(z)\otimes H_j+C_j(z)\otimes K_j.
 \end{aligned}
\end{equation}
The construction proceeds simultaneously with Lemma~\ref{lem: 4.12}.
At each stage that lemma first establishes $\det\Delta_j\not\equiv0$
using the earlier stages. We can then form
$\Delta_j^{-1}\in M_{n_jR_j}(\mathbb C(z))$.

Let $x$ be an admissible tuple on a finite-dimensional Hilbert space
$\mathcal E$, with $d=\dim\mathcal E$. Put
$V_j=I_{n_j}\otimes v_j(x)$ and $B_{P,j}(z)=-A_j(z)$.
Using $E_j(x)=(P_j(x),H_j\otimes I_{\mathcal E},
K_j\otimes I_{\mathcal E})$ with $P_j(x)=v_j(x)v_j(x)^*$, we obtain
\[
 \begin{aligned}
 L_j(z;E_j(x))
 ={}&Z_j(z)\otimes I_{\mathcal K_j}\otimes I_{\mathcal E}
   +A_j(z)\otimes P_j(x)\\
 &+B_j(z)\otimes H_j\otimes I_{\mathcal E}
   +C_j(z)\otimes K_j\otimes I_{\mathcal E}\\
 ={}&\Delta_j(z)\otimes I_{\mathcal E}+A_j(z)\otimes P_j(x)\\
 ={}&\Delta_j(z)\otimes I_{\mathcal E}
   -V_j\bigl(B_{P,j}(z)\otimes I_{\mathcal E}\bigr)V_j^*.
 \end{aligned}
\]
Put
\[
 \mathcal T_j(z;x)=V_j^*(\Delta_j(z)^{-1}\otimes I_{\mathcal E})V_j,
 \qquad
 \alpha_c=I_{n_j}\otimes u_{j,c},\quad
 \beta_c=I_{n_j}\otimes w_{j,c}.
\]
Corollary~\ref{cor 4.5} gives
\[
 \mathcal T_j(z;x)=E_0(z)\otimes I_{\mathcal E}
 +\sum_{r=1}^3
 \bigl(E_r(z)\otimes U(x_r)+E_r^\#(z)\otimes U(x_r)^*\bigr),
\]
where the coefficient matrices in $M_{n_j}(\mathbb C(z))$ are
\[
 \begin{aligned}
 E_0(z)&=\frac16\sum_{c=1}^3
  \left(\alpha_c^*\Delta_j(z)^{-1}\alpha_c+
        \beta_c^*\Delta_j(z)^{-1}\beta_c\right),\\
 E_r(z)&=\frac16\alpha_r^*\Delta_j(z)^{-1}\beta_r
 \qquad (1\leq r\leq3).
 \end{aligned}
\]
Since $\Delta_j^\#=\Delta_j$, we have $E_0^\#=E_0$ and
\[
 E_r^\#(z)=\frac16\beta_r^*\Delta_j(z)^{-1}\alpha_r.
\]
At real $z$ where the coefficients are defined and $\Delta_j(z)$
is invertible, Lemma~\ref{lem: 4.10} gives the equivalence of
invertibility of $L_j(z;E_j(x))$ and
\begin{equation}\label{eq: 4.21}
 F_j(z;x)=
 \begin{pmatrix}
 -\mathcal T_j(z;x)&I_{n_j}\otimes I_{\mathcal E}\\
 I_{n_j}\otimes I_{\mathcal E}&-B_{P,j}(z)\otimes I_{\mathcal E}
 \end{pmatrix}.
\end{equation}
Put $q=2n_j$ and define
\[
 G_0(z)=\begin{pmatrix}-E_0(z)&I_{n_j}\\I_{n_j}&-B_{P,j}(z)\end{pmatrix},
 \qquad
 G_r(z)=\begin{pmatrix}-E_r(z)&0\\0&0\end{pmatrix}
 \quad (1\leq r\leq3).
\]
Then $G_0^\#=G_0$, and $F_j$ is affine in the Cayley unitaries
and their adjoints:
\[
 F_j(z;x)=G_0(z)\otimes I_{\mathcal E}
 +\sum_{r=1}^3
   \bigl(G_r(z)\otimes U(x_r)+G_r^\#(z)\otimes U(x_r)^*\bigr).
\]

If $n_j>\kappa_j$, apply the construction of Lemma~\ref{lem: 4.11}.
Set
\[
 D_0(z)=G_0(z)+\sum_{r=1}^3\bigl(G_r(z)+G_r^\#(z)\bigr),
 \qquad Q_r(z)=\begin{pmatrix}I_q&-2iG_r^\#(z)\end{pmatrix},
\]
and define the formal matrix blocks
\[
 M_r(x)=\begin{pmatrix}
 0&(2-i)I_q+x_rI_q\\
 (2+i)I_q+x_rI_q&0
 \end{pmatrix}.
\]
The next matrix polynomial is
\[
 L_{j+1}(z;x)=
 \begin{pmatrix}
 D_0(z)&Q_1(z)&Q_2(z)&Q_3(z)\\
 Q_1^\#(z)&M_1(x)&0&0\\
 Q_2^\#(z)&0&M_2(x)&0\\
 Q_3^\#(z)&0&0&M_3(x)
 \end{pmatrix},
\]
with scalar blocks amplified by $I_{\mathcal E}$ upon substitution.
At the real points under consideration $Q_r^\#(z)=Q_r(z)^*$,
so Lemma~\ref{lem: 4.11} identifies $F_j(z;x)$ as its Schur
complement with respect to the invertible lower block. Hence
\[
 F_j(z;x)\text{ is invertible}
 \quad\Longleftrightarrow\quad
 L_{j+1}(z;x)\text{ is invertible}.
\]
The new coefficient matrices are formed from the existing
rational matrices and the constants $2,\pm i$. They depend only on
$a$ and the background data at stages $1,\ldots,j$. The block
sizes give $n_{j+1}\leq q+3(2q)=14n_j$.

If $n_j\leq\kappa_j$, Lemma~\ref{lem: 4.3} and
Corollary~\ref{cor 4.5} give $E_r(z)=0$ for $1\leq r\leq3$.
The recursion stops at $J=j$, with $F_J(z)=G_0(z)$ at this stage.
By the size estimate, $n_j\leq14^{j-1}n_1$ at each reached stage.
Since $\kappa_j=j14^{j-1}$, the first stopping stage exists and
$J\leq n_1$. All matrices constructed through stage $J$ depend
only on $a$ and the first $J$ backgrounds.

\medskip
\noindent\emph{Remark.}
Lemma~\ref{lem: 4.8} and Remark~\ref{rem: 4.9} give the initial
reduction. At $x=T_1^{(m)}$ the lower block remains invertible,
and its Schur complement is $zI-a(T_1^{(m)})$. Fix $m\geq J$ and
a real $z$ where the coefficients are defined and the backgrounds
are invertible. The relation $T_j^{(m)}=E_j(T_{j+1}^{(m)})$ gives,
for $j<J$, the chain of invertibility equivalences
\[
 \begin{gathered}
 L_j(z;T_j^{(m)})
 \xleftrightarrow{\text{Lemma~\ref{lem: 4.10}}}
 F_j(z;T_{j+1}^{(m)})
 \xleftrightarrow{\text{Lemma~\ref{lem: 4.11}}}
 L_{j+1}(z;T_{j+1}^{(m)}).
 \end{gathered}
\]
The first step passes to the inner tuple. The second replaces the
dependence on the Cayley unitaries with an affine dependence on
the inner tuple. At stage $J$ the first step still applies, and
the cross compressions vanish. Thus
\[
 F_J\bigl(z;T_{J+1}^{(m)}\bigr)
 =F_J(z)\otimes I_{D_{J+1}^{(m)}}.
\]
The remaining depth determines the size of the identity factor.
\medskip

\begin{lemma}\label{lem: 4.12}
For every stage $j$ reached by the recursion,
$\det\Delta_j(z)$ is a nonzero rational function. The terminal matrix
$F_J(z)\in M_{2n_J}(\mathbb C(z))$ also has nonzero determinant.
\end{lemma}

\begin{proof}
Write $a=\sum_w A_wx_w$, where the sum includes the empty word, and set
$M_a=\sum_w\|A_w\|$. If $S$ is a tuple of contraction matrices,
then $\|S_w\|\leq 1$ for every word $w$; hence
$\|a(S)\|=\|\sum_w A_w\otimes S_w\|\leq M_a.$

For the base case, set $S=(0_{R_1},H_1,K_1)$. By
Remark~\ref{rem: 4.9}, $\Delta_1(z)=L_1(z;S)$ is invertible if and
only if $zI_{\ell R_1}-a(S)$ is invertible. Since $\|a(S)\|\leq M_a$,
this holds for every real $z>M_a$, so $\det\Delta_1(z)\ne0$.

Suppose that $L_j$ has been constructed and that $\det\Delta_i\not\equiv0$ for every $i<j$. Let
$\mathcal B_j\subseteq\mathbb R$ consist of the real poles of the
rational coefficients introduced up to the construction of $L_j$, together with
the real zeros of $\det\Delta_i$ for $i<j$. Since only finitely many coefficients
have been introduced and each earlier determinant is a nonzero
rational function, the set
$\mathcal B_j$ is finite.

Let $S$ be an admissible tuple of $d\times d$ matrices. Set $S_j=S$
and define
\[
S_i=E_i(S_{i+1})~(i=j-1,\ldots,1).
\]
Then $S_1=E_{<j}(S)$ is a tuple of matrices of size
$d_1=d\prod_{i=1}^{j-1}R_i$. For every real $z\notin\mathcal B_j$, the preceding Schur complement equivalences give
\[
 L_i(z;S_i)\text{ is invertible}
 \quad\Longleftrightarrow\quad
 L_{i+1}(z;S_{i+1})\text{ is invertible}
 \qquad(i<j).
\]
Combining these equivalences with Remark~\ref{rem: 4.9} yields
\begin{equation}\label{eq: 4.22}
 \begin{aligned}
 L_j(z;S)\text{ is invertible}
 &\quad\Longleftrightarrow\quad
 L_1(z;S_1)\text{ is invertible}\\
 &\quad\Longleftrightarrow\quad
 zI_{\ell d_1}-a(S_1)\text{ is invertible}.
 \end{aligned}
\end{equation}
For $j=1$, we have $S_1=S$ and $\mathcal B_1=\varnothing$.

To prove that $\det\Delta_j\not\equiv0$, take $S=(0_{R_j},H_j,K_j)$
and choose $z\in(M_a,\infty)\setminus\mathcal B_j$. Since the components of $S_1$
are contractions, it follows that $\|a(S_1)\|\leq M_a<z$ and
$zI_{\ell d_1}-a(S_1)$ is invertible. By \eqref{eq: 4.22}, $\Delta_j(z)=L_j(z; 0,H_j,K_j)$
is invertible for every such $z$. Since $\mathcal B_j$ is finite,
this holds for all sufficiently large real $z$. Hence
$\det\Delta_j\not\equiv0$, and so we may form $\Delta_j^{-1}$ over
$\mathbb C(z)$. If $n_j>\kappa_j$, use it to construct $L_{j+1}$
with $n_{j+1}\leq14n_j\leq14^jn_1$. If $n_j\leq\kappa_j$, the
recursion stops. This completes the induction.

For the terminal determinant, take the coordinate of depth $m=J$.
Its innermost tuple is $T_{J+1}^{(J)}=(0,0,0)\in\mathbb C^3$,
acting on $\mathbb C$, and $T_J^{(J)}=E_J(T_{J+1}^{(J)})$.
Choose a real $z_*>M_a$ outside the finitely many coefficient poles
and background determinant zeros through stage $J$. Since
$\|a(T_1^{(J)})\|\leq M_a$, applying \eqref{eq: 4.22} to
$S=T_J^{(J)}$ gives invertibility of $L_J(z_*;T_J^{(J)})$.
The last Schur complement equivalence then gives invertibility of
\[
 F_J(z_*;T_{J+1}^{(J)})=F_J(z_*)\otimes I_{\mathbb C}.
\]
Thus $\det F_J(z_*)\ne0$, and $\det F_J\not\equiv0$.
\end{proof}

\begin{proposition}\label{prop: 4.13}
For every formally selfadjoint matrix polynomial $a$, there is a
finite set $S_a\subseteq\mathbb R$ such that
$\sigma \left(a(P^{(m)},H^{(m)},K^{(m)})\right)\subseteq S_a
       ~(m\geq1)$.
\end{proposition}

\begin{proof}
Let $J$ be the first stopping stage for $a$. Let $\mathcal E_a$
consist of the real poles of the entries of the scalar rational
matrices used up to stage $J$, including $\Delta_j^{-1}$ and
$F_J$, together with the real zeros of $\det\Delta_j$ for
$1\leq j\leq J$ and of $\det F_J$. By Lemma~\ref{lem: 4.12}, these
determinants are nonzero rational functions. Hence $\mathcal E_a$
is finite and independent of $m$.

Fix $m\geq J$ and $z\in\mathbb R\setminus\mathcal E_a$. For each
$j\leq J$, the tuple $T_{j+1}^{(m)}$ has matrix size
$d_j=\prod_{r=j+1}^{m}R_r$, with the empty product equal to one.
Thus the background matrix at stage $j$ is
$\Delta_j(z)\otimes I_{d_j}$, which is invertible. At the stopping
stage, the terminal matrix is
\[
 F_J(z;T_{J+1}^{(m)})=F_J(z)\otimes I_{d_J},
 \qquad F_J(z)=G_{J,0}(z).
\]

Applying \eqref{eq: 4.22} to $S=T_J^{(m)}$ and then using the
Schur complement equivalence at stage $J$ gives
\[
 \begin{aligned}
 zI_{\ell D_m}-a(T_1^{(m)})\text{ is invertible}
 &\quad\Longleftrightarrow\quad
 L_J(z;T_J^{(m)})\text{ is invertible}\\
 &\quad\Longleftrightarrow\quad
 F_J(z)\otimes I_{d_J}\text{ is invertible}\\
 &\quad\Longleftrightarrow\quad
 F_J(z)\text{ is invertible}.
 \end{aligned}
\]
The last condition holds since $z\notin\mathcal E_a$. Since
$a(T_1^{(m)})$ is selfadjoint, its spectrum is real and therefore
contained in $\mathcal E_a$ for every $m\geq J$.
Only finitely many coordinates with $m<J$ remain, each of finite dimension.
Consequently, the set
\[
 S_a=\mathcal E_a\cup
     \bigcup_{m=1}^{J-1}\sigma \left(a(P^{(m)},H^{(m)},K^{(m)})\right)
\]
is finite and contains the spectrum at every coordinate.
\end{proof}

Using the direct sum in \eqref{eq: 4.17}, we represent $M_q(A)$ on ${\mathcal{H}}_q={\mathbb{C}}^q\otimes\mathcal{H}$.

\begin{proposition}\label{prop: 4.14}
For every $q\geq1$, every selfadjoint element of $M_q(A)$ has
finite spectrum. Moreover, $M_q(A)$ is von Neumann regular.
\end{proposition}
\begin{proof}
Let $a=a^*\in M_q(A)$. Since $A=\mathcal{A}^*(1,p,h,k)$, we may
write $a=F(p,h,k)$ for a matrix polynomial $F$. By replacing $F$ by
$(F+F^*)/2$, we may take $F$ formally selfadjoint and still have $F(p,h,k)=a$.

Write $a_m=F(P^{(m)},H^{(m)},K^{(m)})$ and
$F(x)=\sum_w C_w x_w$, with $C_w\in M_q(\mathbb C)$. Since the coordinate generators are contractions,
\[
\|a_m\|
=\left\|\sum_w C_w\otimes x_w(P^{(m)},H^{(m)},K^{(m)})\right\|
\leq \sum_w\|C_w\|.
\]
For $\xi=(\xi_m)_m$, we have
\[
 \begin{aligned}
 \|a\xi\|^2
 &=\sum_m\|a_m\xi_m\|^2\\
 &\leq\Bigl(\sum_w\|C_w\|\Bigr)^2\sum_m\|\xi_m\|^2
 =\Bigl(\sum_w\|C_w\|\Bigr)^2\|\xi\|^2.
 \end{aligned}
\]
Hence
$a=\bigoplus_ma_m$ is bounded. Proposition~\ref{prop: 4.13} implies there exists
a finite set $S\subseteq\mathbb R$ containing every
$\sigma(a_m)$. Set $Q(t)=\prod_{s\in S}(t-s)$. The spectral
theorem gives $Q(a_m)=0$ at every coordinate, hence $Q(a)=0$.
If $\lambda\notin S$, then
$Q(\lambda)=\prod_{s\in S}(\lambda-s)\neq0$. Since
the polynomial $Q(t)-Q(\lambda)$ vanishes at $t=\lambda$,
there is a polynomial $R_\lambda$ such that
\begin{equation}\label{eq: 4.23}
 \begin{aligned}
 Q(t)-Q(\lambda)=(t-\lambda)R_\lambda(t).
 \end{aligned}
\end{equation}
Evaluating \eqref{eq: 4.23} at $a$ gives $(a-\lambda I)R_\lambda(a)
 =Q(a)-Q(\lambda)I
 =-Q(\lambda)I.$
Since $R_\lambda(a)$ commutes with $a$, we have $R_\lambda(a)(a-\lambda I)=-Q(\lambda)I$
as well. Set
$B_\lambda=-\frac{R_\lambda(a)}{Q(\lambda)}.$
The operator $B_\lambda$ is bounded because it is a polynomial in $a$.
Therefore
\[
 (a-\lambda I)B_\lambda=B_\lambda(a-\lambda I)=I.
\]
Thus $a-\lambda I$ has a bounded two-sided inverse. Since this holds
for every $\lambda\in\mathbb C\setminus S$, we conclude that
$\sigma(a)\subseteq S$.

For arbitrary $b\in M_q(A)$, apply the first assertion to
$x=b^*b$. Interpolate a polynomial $f(t)$ on $\sigma(x)\cup\{0\}$, with $f(0)=0$ and $f(t)=t^{-1}$ at every nonzero
spectral value. Put $e=xf(x)$ and $c=f(x)b^*$.
The operator $e$ is the spectral projection onto $(\ker x)^\perp$,
and $x(1-e)=0$. Hence
$\|b(1-e)\xi\|^2=\langle x(1-e)\xi,(1-e)\xi\rangle=0$,
so $be=b$. It follows that
\[
 bcb=bf(x)b^*b=bf(x)x=be=b.
\]
Thus $M_q(A)$ is regular.
\end{proof}

\begin{lemma}\label{lem: 4.15}
Let $D$ be a subalgebra of $B(\mathcal E)$ containing the identity $I_{\mathcal E}$. If $D$ is
von Neumann regular and has at most countable dimension over $\mathbb{C}$,
then every element of $D$ is algebraic.
\end{lemma}
\begin{proof}
Fix $b\in D$. If $y=b-\lambda I$ is invertible in $B(\mathcal E)$,
regularity gives $d\in D$ with $ydy=y$. Multiplication on both
sides by $y^{-1}$ gives $d=y^{-1}$. In particular,
$(b-\lambda I)^{-1}\in D$ for every real $\lambda>\|b\|$.

Suppose that $b$ is not algebraic. We claim that these uncountably
many resolvents are linearly independent. If $\lambda_1,\ldots,\lambda_m$
are distinct and $\sum_j\alpha_j(b-\lambda_jI)^{-1}=0$,
multiply by $\prod_j(b-\lambda_jI)$ to obtain $P(b)=0$, where
\[
 P(t)=\sum_{j=1}^m\alpha_j\prod_{i\ne j}(t-\lambda_i).
\]
The assumption on $b$ implies $P\equiv0$. Evaluating at
$t=\lambda_j$ gives $\alpha_j\prod_{i\ne j}(\lambda_j-\lambda_i)=0$,
so $\alpha_j=0$ for every $j$. Thus the resolvents are linearly
independent, contradicting the countable complex dimension of $D$.
Therefore $b$ is algebraic.
\end{proof}

\begin{proposition}\label{prop: 4.16}
The algebra $A$ is an infinite-dimensional $R^*$-algebra. For every
$q\geq1$, every element of $M_q(A)$ is algebraic over $\mathbb{C}$.
\end{proposition}

\begin{proof}
The coordinate maps are surjective onto $M_{D_m}(\mathbb C)$
by Lemma~\ref{lem: 4.6}. Since $D_m\geq60^m$, the algebra $A$ is
infinite-dimensional.
Proposition~\ref{prop: 4.14} shows that $A$ and $M_q(A)$
are regular. Finite words in $p,h,k$ span $A$, so $M_q(A)$
has at most countable complex dimension. Apply Lemma~\ref{lem: 4.15}
to $M_q(A)$ on $\mathcal H_q$.
\end{proof}

\begin{example}\label{ex:three-projections}
Set $x_1=p$, $x_2=h$, $x_3=k$, and $U_c=U(x_c)$.
The Lagrange interpolation formula on the finite spectrum of each
$x_c$ gives $U_c\in A$ by \eqref{eq: 4.10}. The recovery formula
\eqref{eq: 4.11} for the Cayley unitaries then gives
\[
 A=\mathcal A^*(1,p,h,k)=\mathcal A^*(1,U_1,U_2,U_3).
\]
Let $t_1=1/4$, $t_2=1/2$, $t_3=3/4$, and
$\alpha_c=\sqrt{t_c(1-t_c)}$. On ${\mathcal H}'={\mathbb C}^6\otimes\mathcal{H}$ define
\begin{equation}\label{eq:three-projections}
 \begin{aligned}
 \eta&=\frac1{\sqrt6}(1,U_1,1,U_2,1,U_3)^T,
       &\mathsf P&=\eta\eta^*,\\
 \mathsf Q&=\bigoplus_{c=1}^3\begin{pmatrix}1&0\\0&0\end{pmatrix},
       &\mathsf R&=\bigoplus_{c=1}^3
            \begin{pmatrix}1-t_c&\alpha_c\\\alpha _c&t_c\end{pmatrix}.
 \end{aligned}
\end{equation}
Scalar entries are multiples of $I_{\mathcal H}$. The operator
$\mathsf Q$ is a projection. If $R_c$ is the $c$th block of
$\mathsf R$, then $\alpha_c^2=t_c(1-t_c)$ gives
\[
 R_c^2=\begin{pmatrix}
 (1-t_c)^2+\alpha_c^2&\alpha_c(1-t_c)+\alpha_ct_c\\
 \alpha_c(1-t_c)+\alpha_ct_c&\alpha_c^2+t_c^2
 \end{pmatrix}=R_c,
\]
Thus $\mathsf R$ is a projection. Also
$\eta^*\eta=6^{-1}\sum_c(1+U_c^*U_c)=1$, so $\mathsf P$ is a projection.
Its matrix is
\[
\mathsf P=\frac16
\begin{pmatrix}
I   & U_1^*      & I   & U_2^*      & I   & U_3^*\\
U_1 & I          & U_1 & U_1U_2^*  & U_1 & U_1U_3^*\\
I   & U_1^*      & I   & U_2^*      & I   & U_3^*\\
U_2 & U_2U_1^*   & U_2 & I          & U_2 & U_2U_3^*\\
I   & U_1^*      & I   & U_2^*      & I   & U_3^*\\
U_3 & U_3U_1^*   & U_3 & U_3U_2^*  & U_3 & I
\end{pmatrix}.
\]

Let $D=\mathcal A^*(\mathsf P,\mathsf Q,\mathsf R)\subseteq M_6(A)$
and set $S=(\mathsf Q-\mathsf R)^2$. Its $c$th block is
\[
 \begin{pmatrix}t_c&-\alpha_c\\-\alpha_c&-t_c\end{pmatrix}^{\!2}
 =t_c I_2,
 \qquad S=\bigoplus_{c=1}^3 t_c I_2.
\]
The Lagrange interpolation formula at $0,t_1,t_2,t_3$ gives
\[
 f_c(t)=\frac{t}{t_c}\prod_{d\ne c}\frac{t-t_d}{t_c-t_d},
 \qquad f_c(0)=0,\quad f_c(t_d)=\delta_{cd}.
\]
Thus the projection $F_c=f_c(S)$ onto the $c$th block belongs to
$D$, and $\sum_c F_c=I_6\otimes I_{\mathcal H}$.
Put $o_c=2c-1$, $e_c=2c$, and write
$\widehat E_{ij}=E_{ij}\otimes I_{\mathcal H}$.
The diagonal matrix units in the $c$th block belong to $D$, since
\[
 \widehat E_{o_co_c}=F_c\mathsf QF_c,\qquad
 \widehat E_{e_ce_c}=F_c-\widehat E_{o_co_c}=F_c-F_c\mathsf QF_c.
\]
Moreover,
\[
 \widehat E_{o_co_c}\mathsf R\widehat E_{e_ce_c}
 =\alpha_c\widehat E_{o_ce_c},
 \qquad
 \widehat E_{o_co_c}\mathsf P\widehat E_{o_do_d}
 =\frac16\widehat E_{o_co_d}.
\]
Since $\alpha_c\ne0$, these identities and adjoints give
\[
 \widehat E_{ij}
 =\widehat E_{i,o_c}\widehat E_{o_c,o_d}\widehat E_{o_d,j}\in D
 \quad\bigl(i\in\{o_c,e_c\},\ j\in\{o_d,e_d\}\bigr).
\]

The $(1,e_d)$ entry of $\mathsf P$ is $U_d^*/6$. Thus
$6\widehat E_{11}\mathsf P\widehat E_{e_d,1}=E_{11}\otimes U_d^*\in D$.
Taking adjoints gives $E_{11}\otimes U_d\in D$. These elements
and the corner identity $\widehat E_{11}$ generate
$E_{11}\otimes A$. For $a\in A$ we therefore have
\[
 \widehat E_{i1}(E_{11}\otimes a)\widehat E_{1j}
 =E_{ij}\otimes a\in D.
\]
Their finite sums fill $M_6(A)$, so $D=M_6(A)$.

\end{example}

Summarizing the above example, we obtain the following theorem.
\begin{theorem}\label{thm:main}
There are three projections $\mathsf P,\mathsf Q,\mathsf R$ on a
complex Hilbert space such that the algebraic $*$-algebra
$\mathcal B=\mathcal{A}^*(\mathsf P,\mathsf Q,\mathsf R)$ is unital and
infinite-dimensional. For every $q\geq1$, $M_q(\mathcal B)$ is
an $R^*$-algebra, and every element of $M_q(\mathcal B)$ is
algebraic over $\mathbb{C}$.
\end{theorem}
\begin{proof}
Example~\ref{ex:three-projections} identifies $\mathcal B$ with
$M_6(A)$ and supplies its unit and three generators. It is
infinite-dimensional by Proposition~\ref{prop: 4.16}.
For each $q\geq1$, the canonical matrix identification
$M_q(M_6(A))\cong M_{6q}(A)$ preserves multiplication and adjoints.
Propositions~\ref{prop: 4.14} and \ref{prop: 4.16} give regularity
and algebraicity at this matrix level.
\end{proof}

\begin{corollary}\label{cor:II1-embedding}
The three-projection example admits a unital $*$-embedding into
every type $\mathrm{II}_1$ factor.
\end{corollary}
\begin{proof}
The representation identifies $\mathcal B$ with a unital
$*$-subalgebra of the bounded von Neumann algebra product
\[
 \mathcal N=\prod_{m\geq1}M_{6D_m}(\mathbb C).
\]
Fix a type $\mathrm{II}_1$ factor $\mathcal M$ with normalized trace $\tau$. By the continuous
range of $\tau$ on projections, choose pairwise orthogonal projections
$e_m\in\mathcal M$ with $\tau(e_m)=2^{-m}$. Normality and faithfulness
of $\tau$ give $\sum_m e_m=1$ in the strong operator topology.

Fix $m$ and put $e=e_m$ and $r=6D_m$. In $e\mathcal M e$,
split $e$ into $r$ orthogonal projections $q_a$ of trace $\tau(e)/r$.
Since equal-trace projections in a factor are equivalent, choose
$v_a$ with $v_a^*v_a=q_1$ and $v_av_a^*=q_a$. The operators
$f_{ab}=v_av_b^*$ satisfy
\[
 f_{ab}f_{dj}=\delta_{bd}f_{aj},\qquad
 f_{ab}^*=f_{ba},\qquad
 \sum_{a=1}^r f_{aa}=e_m.
\]
They give the matrix-unit identification
\[
\iota_m\colon M_r(\mathbb C)
\xrightarrow{\ \cong\ }
\operatorname{span}_{\mathbb C}\{f_{ab}:1\leq a,b\leq r\}
\subseteq e\mathcal M e,\qquad
[a_{ab}]\longmapsto\sum_{a,b=1}^r a_{ab}f_{ab}.
\]

Compression by the $e_m$ identifies the block-diagonal von Neumann
algebra with the product of the corners:
\[
 \mathcal M\cap\{e_m:m\geq1\}'
 \ \cong\ \prod_{m\geq1}e_m\mathcal M e_m,
 \qquad x\longmapsto(e_mxe_m)_m.
\]
For a bounded family $(x_m)_m$ in the corners and $\xi$ in the
representation space,
\[
 \left\|\sum_{m=N'+1}^{N''}x_m\xi\right\|^2
 \leq\left(\sup_m\|x_m\|\right)^2
      \sum_{m=N'+1}^{N''}\|e_m\xi\|^2.
\]
The right side tends to zero as $N',N''\to\infty$, and the strong
limit lies in $\mathcal M$. Its compressions are the given $x_m$,
which gives the inverse identification. Composing it with the product of the
$\iota_m$ gives a normal unital $*$-embedding
\[
 \Phi\colon\mathcal N\longrightarrow\mathcal M,\qquad
 \Phi((a_m)_m)=\operatorname*{s\!-\!lim}_{N\to\infty}
                  \sum_{m=1}^{N}\iota_m(a_m).
\]
Restricting $\Phi$ to $\mathcal B$, we obtain
\[
 \Phi(\mathcal B)
 =\operatorname{alg}^*\bigl(
      \Phi(\mathsf P),\Phi(\mathsf Q),\Phi(\mathsf R)\bigr)
 \subseteq\mathcal M.
\]
The three images are projections, and $\Phi(\mathcal B)$ is an
infinite-dimensional $R^*$-algebra because $\Phi|_{\mathcal B}$
is a $*$-isomorphism onto its image.
\end{proof}

\section*{Acknowledgments}
The development and refinement of the mathematical arguments in this
manuscript benefited from discussions with ChatGPT (GPT-6). All AI-assisted contributions
have been carefully reviewed and verified by the authors, who take full responsibility
for the content of the manuscript.

\section*{Conflict of interest}
The authors declare no conflict of interest.

\end{document}